\documentclass[10pt]{amsart}
\usepackage[british,UKenglish]{babel}
\usepackage[utf8]{inputenc}
\usepackage{braket}
\usepackage{xfrac}
\usepackage{amsfonts, amssymb, amsmath, amsthm, amscd}
\usepackage{mathtools}
\usepackage{tikz}
\usepackage{empheq}
\usepackage[left=3.3cm,right=3.3cm,top=3.3cm,bottom=3.3cm]{geometry}
\usepackage{color}
\usepackage{epsfig}  	
\usepackage{epic,eepic}   
\usepackage{setspace}
\usepackage{caption}
\usepackage{graphicx}
\newcommand{\norm}[1]{\left\lVert #1 \right\rVert}
\newtheorem{thm}{Theorem}[section]
\newtheorem{thmA}{Theorem}
\newtheorem{corA}[thmA]{Corollary}

\newtheorem{prop}[thm]{Proposition}
\newtheorem{lem}[thm]{Lemma}

\theoremstyle{definition}

\theoremstyle{remark}
\newtheorem{remark}[thm]{Remark}

\DeclareMathOperator{\Vol}{Vol}

\DeclareMathOperator{\Div}{div}

\DeclareMathOperator{\interior}{int}
\DeclareMathOperator{\id}{id}
\newcommand{\dmu}{\, \mathrm{d}\mu}
\newcommand{\dV}{\, \mathrm{d}V}
\newcommand{\dA}{\, \mathrm{d}A}
\newcommand{\ddt}{\frac{\mathrm{d}}{\mathrm{d}t}}
\newcommand{\dd}{\mathrm{d}}

\newcommand{\Span}{\mathrm{span}}
\begin{document}

\makeatletter
\@namedef{subjclassname@2020}{%
\textup{2020} Mathematics Subject Classification}
\makeatother

\begin{abstract}
We study the $G$-invariant Steklov spectrum of compact connected Riemannian manifolds with boundary, endowed with an effective isometric action of a compact Lie group $G$. Under a mild nontransitivity assumption on the boundary action, we establish that this spectrum consists of an infinite sequence of eigenvalues diverging to infinity. We then investigate extremal metrics for the normalised $G$-invariant Steklov eigenvalue functionals. Assuming that $\dim G \geq 1$, we prove that none of these functionals admits an extremal $G$-invariant metric. Moreover, when $\dim M \geq 3$, the normalised first $G$-invariant Steklov eigenvalue is unbounded above even within a fixed $G$-invariant conformal class. The proof relies on a local model near a boundary orbit and spectral convergence under conformal collapse. Finally, in contrast with the upper bounds obtained by Colbois, Dryden and El Soufi for invariant Laplace eigenvalues of embedded spheres of revolution, we construct, for any $n \geq 3$, a sequence of smooth Euclidean domains of revolution in $\mathbb R^n$, all diffeomorphic to the ball, along which the normalised first $G$-invariant Steklov eigenvalue blows up.
\end{abstract}

\title[Extremal metrics for the $G$-invariant Steklov eigenvalues]{Extremal metrics for the $G$-invariant \\ Steklov eigenvalues}
\author{Eduardo Longa}

\address{Departamento de Matem\'{a}tica, Instituto de Matem\'{a}tica, Estat\'{i}stica e Ci\^encia da Computa\c{c}\~ao, Universidade de S\~{a}o Paulo, R. do Mat\~{a}o 1010, S\~{a}o Paulo, SP 05508-900, Brazil}
\email{edulonga@ime.usp.br}



\maketitle

\section{Introduction} \label{sec: Introduction}

Let $(M, g)$ be a compact connected Riemannian manifold of dimension $n \geq 2$ with smooth nonempty boundary. The Dirichlet-to-Neumann operator $\mathcal{D}_g : H^{1/2}(\partial M) \to H^{-1/2}(\partial M)$ is defined by 
\begin{align*}
\mathcal{D}_g(f) = \frac{\partial \hat f}{\partial \nu},
\end{align*}
where $\hat f \in H^1(M)$ is the (weakly) harmonic extension of $f$ to $M$ and $\nu$ is the outward unit normal to $\partial M$. The eigenvalues of $\mathcal{D}_g$ are called the Steklov eigenvalues of $(M,g)$, and they form a diverging sequence of nonnegative numbers
\begin{align*}
0 = \sigma_0(g) < \sigma_1(g) \leq \cdots  \leq \sigma_k(g) \leq \cdots \nearrow + \infty,
\end{align*} 
where each eigenvalue is repeated according to its multiplicity. The associated eigenfunctions $f_k \in H^{1/2}(\partial M)$ are known to be smooth by standard elliptic theory, and they form an orthonormal basis of $L^2(\partial M)$ after suitable normalisation. Their harmonic extensions $u_k = \hat f_k \in C^\infty(M)$ are solutions to the Steklov eigenvalue problem
\begin{align*}
\begin{cases}
\Delta_g u_k = 0 & \text{in } M,\\
\partial_\nu u_k = \sigma_k(g) u_k & \text{on } \partial M,
\end{cases}
\end{align*}
where $\Delta_g = -\Div_g \nabla$ is the Laplace--Beltrami operator. We refer to the functions $u_k$ as Steklov eigenfunctions, although some authors prefer to reserve this term for the boundary traces $f_k$. 

The search for sharp eigenvalue bounds under geometric constraints is a classical theme in spectral geometry.
For the Steklov problem, this line of research goes back at least to the celebrated work of Weinstock \cite{weinstock} in 1954, who proved that every bounded simply connected planar domain $\Omega$ satisfies $\sigma_1(\Omega) L(\partial \Omega) \leq 2 \pi$, with equality if and only if $\Omega$ is a disc; here, $L(\partial \Omega)$ denotes the boundary length. In 1974, Hersch--Payne--Schiffer \cite{hersch} extended this estimate to higher eigenvalues, showing that $\sigma_k(\Omega) L(\partial \Omega) \leq 2 k \pi$ for every such domain and every $k \geq 1$. Girouard--Polterovich \cite{girouard2010} subsequently proved that this bound is sharp, exhibiting sequences of simply connected domains degenerating to a disjoint union of $k$ identical discs along which the upper bound is approached. More recently, Girouard--Karpukhin--Lagacé \cite{girouard2021} solved the corresponding optimisation problem without the assumption of simple connectivity, proving the sharp inequality $\sigma_k(\Omega) L(\partial \Omega) \leq 8 k \pi$ for bounded planar domains.

These optimisation problems have also been studied in the broader setting of compact Riemannian surfaces with boundary. In this direction, Fraser--Schoen \cite{fraser2011} generalised Weinstock's inequality, proving that $\sigma_1(g) L(\partial \Sigma) \leq 2 (\gamma + b) \pi$ for every compact connected orientable Riemannian surface $(\Sigma, g)$ of genus $\gamma$ with $b$ boundary components. This result was later extended to higher eigenvalues by Girouard--Polterovich \cite{girouard2012} and subsequently improved by Karpukhin \cite{karpukhin2017}, who showed that $\sigma_k(g) L(\partial \Sigma) \leq 2 \pi (\gamma + b + k - 1)$ under the same hypotheses.
 
 Another fundamental problem is to investigate the existence and characterisation of extremal metrics for Steklov eigenvalues. A remarkable result of Fraser--Schoen \cite{fraser} establishes a connection with minimal surface theory: metrics maximising normalised Steklov eigenvalues on compact Riemannian surfaces with boundary arise from free boundary minimal immersions into Euclidean unit balls. This can be viewed as a Steklov analogue of Nadirashvili's theorem \cite{nadirashvili}, which relates maximising metrics for the first area-normalised Laplace eigenvalue on closed surfaces to minimal immersions into Euclidean spheres.
 
Recently, Karpukhin--M\'etras \cite{karpukhin2022} extended these geometric characterisations to arbitrary dimensions and to extremal metrics that are not necessarily maximisers, in terms of free boundary minimal immersions and $n$-harmonic maps. Their results are formulated using the normalisation
 \begin{align*}
 \overline{\sigma}_k(g) = \sigma_k(g) \Vol(M, g)^{\frac{2 - n}{n}} \Vol(\partial M, g).
 \end{align*}
This choice is natural for two reasons. First, in dimension two, it reduces to the usual boundary length normalisation. Second, in every dimension $n \geq 2$, it is the only scale-invariant normalisation involving powers of $\Vol(M, g)$ and $\Vol(\partial M, g)$ that allows the existence of extremal metrics \cite[Proposition 1]{karpukhin2022}.

Many other works have investigated the optimisation of Steklov eigenvalues and related questions in spectral geometry. For a broader discussion of these and other recent developments, we refer to the surveys \cite{colbois survey, girouard}. For a historical overview of the Steklov problem and its connections to mathematical physics, we refer to \cite{kuznetsov}.

A natural question is how these problems in spectral geometry change when symmetry is imposed not only on the metrics, but also on the eigenfunctions. Suppose that a compact Lie group $G$ acts smoothly on $M^n$. For each $G$-invariant Riemannian metric $g$, the Dirichlet-to-Neumann operator commutes with the $G$-action. Restricting the Steklov problem to $G$-invariant functions defines the $G$-invariant Steklov eigenvalues, denoted by $\sigma_k^G(g)$. These form a subsequence of the ordinary Steklov spectrum, with each eigenvalue repeated according to the dimension of the space of associated $G$-invariant eigenfunctions. Under a mild nontransitivity assumption on the boundary action, specified below, these eigenvalues constitute an infinite diverging sequence
\begin{align*}
0 = \sigma_0^G(g) < \sigma_1^G(g) \leq \sigma_2^G(g) \leq \cdots \leq \sigma_k^G(g) \leq \cdots \nearrow \infty.
\end{align*}
In particular, $\sigma_k(g) \leq \sigma_k^G(g)$ for every $k \geq 0$, although the position of $\sigma_k^G(g)$ in
the ordinary spectrum may depend on $g$. With the same normalisation as above, we define
\begin{align*}
\overline{\sigma}_k^G(g) = \sigma_k^G(g) \Vol(M, g)^{\frac{2 - n}{n}}  \Vol(\partial M, g).
\end{align*}
These are the functionals studied in this paper. Their optimisation differs from that of the ordinary Steklov eigenvalues among $G$-invariant metrics, since the eigenfunctions must also be $G$-invariant. On closed manifolds, the analogous restriction of the Laplace--Beltrami operator allows one to define the $G$-invariant Laplace eigenvalues, denoted by $\lambda_k^G(g)$, and the corresponding normalised eigenvalues $\overline{\lambda}_k^G(g) = \lambda_k^G(g) \Vol(M, g)^{2/n}$.

Early work on the optimisation of invariant Laplace eigenvalues was carried out by Engman \cite{engman} and Abreu--Freitas \cite{abreu}, in the setting of $\mathbb{S}^1$-invariant metrics on the two-dimensional sphere $\mathbb{S}^2$. Engman proved that the first positive $\mathbb{S}^1$-invariant eigenvalue is unbounded above among metrics of fixed area, while Abreu--Freitas showed that $\lambda_1^{\mathbb{S}^1}(g)$ can in fact attain any positive value when restricted to metrics of fixed area. They also established sharp upper bounds for all invariant eigenvalues when the metrics are induced by embedded surfaces of revolution in $\mathbb{R}^3$.
In particular, for such metrics of area $4 \pi$, $\lambda_1^{\mathbb{S}^1}(g) < \frac{\xi_1^2}{2}$, where $\xi_1$ is the first positive zero of the Bessel function $J_0$, and the constant on the right-hand side is optimal. 

Colbois--Dryden--El Soufi \cite{colbois_G} subsequently extended this study to a general closed connected manifold $M^n$ endowed with an effective nontransitive action of a compact Lie group $G$. They proved that if $\dim G \geq 1$, the functionals $\overline{\lambda}_k^G$ admit no extremal $G$-invariant metrics. Moreover, assuming $n \geq 3$, the normalised first $G$-invariant eigenvalue $\overline{\lambda}_1^G$ is unbounded above even within a fixed $G$-invariant conformal class. This contrasts with Korevaar's \cite{korevaar} boundedness result for the ordinary Laplace eigenvalues in a conformal class. For the standard $\mathrm{O}(n)$-action on $\mathbb{S}^n$, they also established sharp upper bounds for $\overline{\lambda}_k^{\mathrm{O}(n)}$ when the metrics are induced by embeddings as hypersurfaces of revolution in $\mathbb{R}^{n+1}$.

In this paper, we investigate the corresponding questions for the $G$-invariant Steklov spectrum. Unlike in the closed setting, nontransitivity of the action on $M$ does not by itself guarantee that the invariant spectrum is infinite. For instance, let 
\begin{align*}
(M,g) = \left( [0,L] \times \mathbb{S}^{n - 1}, \dd t^2 \oplus g_{\mathbb{S}^{n - 1}} \right),
\end{align*}
and let $G = \mathrm{SO}(n)$ act trivially on the first factor and by the standard action on $\mathbb{S}^{n - 1}$. Since the action of $\mathrm{SO}(n)$ on $\mathbb{S}^{n - 1}$ is transitive, every $G$-invariant function on $M$ depends only on the variable $t$. Thus, every $G$-invariant harmonic function is of the form $u(t, x) = a + bt$, and a direct computation shows that the $G$-invariant Steklov spectrum consists only of $\{ 0, \frac{2}{L} \}$. 

To exclude this phenomenon, we impose the following condition. For each connected component $N$ of $\partial M$, let
\begin{align*}
G_N = \{ \gamma \in G : \gamma(N) = N \}
\end{align*}
be the subgroup preserving $N$. We call the action \textbf{boundary-nontransitive} if the induced action of $G_N$
on $N$ is nontransitive for at least one connected component $N$ of $\partial M$. As shown in Section \ref{sec: Steklov and Steklov-Neumann eigenvalues}, this condition is equivalent to the existence of infinitely many linearly
independent smooth $G$-invariant functions on $\partial M$ and ensures that the $G$-invariant Steklov spectrum is infinite. When $\partial M$ is connected, this condition is equivalent to nontransitivity of the $G$-action on $\partial M$. For disconnected boundaries, it is stronger: there must exist two points in the same connected component of $\partial M$ that belong to distinct $G$-orbits.

We first establish the nonexistence of extremal metrics for the normalised $G$-invariant Steklov eigenvalues.
Here extremality is understood with respect to smooth deformations through $G$-invariant metrics; the precise
definition is given in Section \ref{sec: Extremal metrics for the Steklov eigenvalues}.

\begin{thmA} \label{thm: no extremal metrics}
Let $M$ be a compact connected smooth manifold with boundary $\partial M \neq \emptyset$ and $\dim M = n \geq 3$. Suppose $G$ is a compact Lie group acting smoothly and effectively on $M$, and assume the action is boundary-nontransitive. If $\dim G \geq 1$, then there are no $\overline{\sigma}_k^G$-extremal $G$-invariant Riemannian metrics on $M$, for any $k \geq 1$.
\end{thmA}

\begin{remark} \label{remark: n geq 3 in main theorem}
Under the other hypotheses of Theorem \ref{thm: no extremal metrics}, the assumption $n \geq 3$ is automatic.	Indeed, if a compact Lie group of positive dimension acts smoothly and effectively on a compact connected surface	with boundary, then $G_N$ acts transitively on all boundary components $N \subset \partial M$. Thus, such an action cannot be boundary-nontransitive. We prove this in Remark \ref{remark: n geq 3 proof}.
\end{remark}

Our second main result provides a Steklov analogue of the unboundedness theorem of Colbois--Dryden--El Soufi
\cite[Theorem~1.2]{colbois_G}, showing that $\overline{\sigma}_1^G$ is unbounded above even within a fixed conformal class.

\begin{thmA} \label{thm: infinite sup}
Let $(M, g_0)$ be a compact connected smooth Riemannian manifold with boundary $\partial M \neq \emptyset$ and $\dim M = n \geq 3$. Suppose $G$ is a compact Lie group of dimension $\dim G \geq 1$ acting smoothly and effectively on $(M, g_0)$ by isometries, and assume the action is boundary-nontransitive. Then
\begin{align*}
\sup_g \sigma_1^G(g) \Vol(M, g)^{\frac{2 - n}{n}} \Vol(\partial M, g) = + \infty,
\end{align*}
where the supremum is taken over all smooth metrics $g$ that are $G$-invariant and conformal to $g_0$.
\end{thmA}

This contrasts with the ordinary Steklov spectrum. Hassannezhad \cite[Theorem 4.1]{hassannezhad} proved
that, for each fixed conformal class $[g_0]$ of metrics on $M$, there exists a constant $C > 0$ such that
\begin{align*}
\overline{\sigma}_k(g) \leq C k^{2/n}
\end{align*}
for every $g \in [g_0]$ and every $k \geq 1$. Thus, although the ordinary normalised Steklov eigenvalues
remain bounded in $[g_0]$ for each fixed $k$, the first positive normalised $G$-invariant Steklov eigenvalue can become arbitrarily large even among $G$-invariant metrics in this conformal class. As in \cite[Corollary 1.4]{colbois_G}, combining this estimate with Theorem~\ref{thm: infinite sup} yields the following consequence.

\begin{corA} \label{cor: first eigenfunctions}
Let $(M,g_0)$ and $G$ be as in Theorem \ref{thm: infinite sup}. For every positive integer $K$, there exists a smooth $G$-invariant metric $g_K$ conformal to $g_0$ such that none of the first $K$ nonconstant Steklov eigenfunctions of $g_K$ is $G$-invariant.
\end{corA}

Phrased differently, it is possible to conformally deform $g_0$ among $G$-invariant metrics so that the first positive $G$-invariant Steklov eigenvalue appears at an arbitrarily large index in the ordinary Steklov spectrum.

Finally, we investigate whether imposing an extrinsic Euclidean constraint changes the unboundedness phenomenon of Theorem~\ref{thm: infinite sup}. The upper bounds for invariant Laplace eigenvalues of embedded hypersurfaces of revolution discussed above naturally raise the question of whether analogous bounds hold for the invariant Steklov spectrum of Euclidean domains of revolution. In this setting, we keep the ambient Euclidean metric fixed and vary the domain.

We consider solid Euclidean domains of revolution, rather than hypersurfaces of revolution with boundary. This distinction is important: the rotation group acts transitively on the boundary of a rotational hypersurface diffeomorphic to a ball, whereas its action on the boundary of a solid domain of revolution is nontransitive.

\begin{thmA} \label{thm: revolution balls unbounded}
Let $n \geq 3$ and let $G = \mathrm{O}(n - 1)$ act on $\mathbb{R}^n = \mathbb{R}^{n - 1} \times \mathbb{R}$ by orthogonal transformations on the first factor. Then
\begin{align*}
\sup_B \sigma_1^G(B) \Vol(B)^{\frac{2 - n}{n}} \Vol(\partial B) = +\infty,
\end{align*}
where the supremum is taken over all smooth bounded $G$-invariant domains $B \subset \mathbb{R}^n$ diffeomorphic to the ball, equipped with the Euclidean metric.
\end{thmA}
The proof is presented in Section~\ref{sec: Revolution-domains}, where we first treat the three-dimensional case in detail, as the construction is easier to visualise, and then explain the extension to higher dimensions.

We briefly describe the main ideas behind our proofs, which build on the techniques introduced by Colbois--Dryden--El Soufi \cite{colbois_G} in the closed setting. The nonexistence result follows from a variational characterisation of extremal metrics presented in Section~\ref{sec: Characterisation of extremal metrics and proof of Theorem A}, which implies that the metric can be expressed as a sum of squares of differentials of $G$-invariant Steklov eigenfunctions. Since these functions are constant along $G$-orbits, their differentials vanish in tangent directions, contradicting the positive definiteness of the metric whenever the action has positive-dimensional orbits. To establish unboundedness of $\overline{\sigma}_1^G$ within a conformal class, we construct shrinking neighbourhoods of boundary orbits on which the normalised $G$-invariant Steklov--Neumann eigenvalues diverge, and then transfer this behaviour to the Steklov spectrum of $M$ through conformal collapse. Finally, the Euclidean result is obtained by constructing domains of revolution with rapidly oscillating boundaries, whose boundary areas grow sufficiently fast to force $\overline{\sigma}_1^G$ to diverge, despite the decay of the eigenvalue itself.

The paper is organised as follows. In Section \ref{sec: Steklov and Steklov-Neumann eigenvalues}, we establish the spectral preliminaries for the $G$-invariant Steklov and Steklov--Neumann problems. Section \ref{sec: Extremal metrics for the Steklov eigenvalues} develops the variational framework for extremal metrics, while Section \ref{sec: Characterisation of extremal metrics and proof of Theorem A} gives a characterisation of $G$-invariant extremal metrics and proves Theorem \ref{thm: no extremal metrics}. In Section \ref{sec: The proof of Theorem infinite sup}, we prove Theorem \ref{thm: infinite sup} using local models near boundary orbits and conformal collapse, and deduce Corollary C. The spectral
convergence results used in the proof are established in Section  \ref{sec: conformal localisation}. Finally, Section \ref{sec: Revolution-domains} is devoted to the construction of Euclidean domains of revolution with unbounded $\overline{\sigma}_1^G$.

\bigskip

\textbf{Statement on the use of AI.} OpenAI's ChatGPT (GPT-5.6 Sol and GPT-6 Astra) was used during the preparation of this manuscript for grammar checking, overall organisation, proofreading, and figure preparation. The estimates in the proof of Proposition~5.4 were developed through interaction with the tool. It also proposed the initial construction of the Euclidean domains of revolution in Section \ref{sec: Revolution-domains}, which was subsequently refined and analysed through the same process. The author independently checked the mathematical arguments and calculations in these parts of the manuscript and takes full responsibility for the content of the manuscript.

\section*{Acknowledgements}

The author expresses his gratitude to Bruno Colbois for kindly suggesting the main idea behind the proof of Theorem~\ref{thm: infinite sup}. He also thanks Pieralberto Sicbaldi for stimulating discussions and for his hospitality during a visit to the University of Granada. The author was partially supported by grant 2024/01663-8, São Paulo Research Foundation (FAPESP).

\section{Steklov and Steklov--Neumann eigenvalues} \label{sec: Steklov and Steklov-Neumann eigenvalues}

In this section, we collect the spectral preliminaries that will be used throughout the paper. We first recall the Steklov problem and introduce its $G$-invariant counterpart, including a criterion ensuring that the $G$-invariant spectrum is infinite. We then consider the mixed Steklov--Neumann problem and its $G$-invariant version, which will play a central role in the proof of Theorem \ref{thm: infinite sup}.

\subsection{The Steklov problem} \label{sec: The Steklov problem}

We begin by recalling the standard Steklov problem and its variational formulation. We then restrict the problem to $G$-invariant functions when a compact Lie group acts by isometries on the manifold.

Let $(M, g)$ be a compact connected smooth Riemannian manifold with nonempty boundary and dimension $n \geq 2$, and denote by $\Delta_g = - \Div_g \nabla$ its Laplace--Beltrami operator.

We denote by $H^{1/2}(\partial M)$ the Sobolev space defined using the Laplace--Beltrami operator $\Delta_{\partial M}$ of the metric induced by $g$ on $\partial M$. More precisely, if $\{\varphi_j\}_{j=0}^\infty$ is an orthonormal basis of $L^2(\partial M)$ satisfying $\Delta_{\partial M}\varphi_j = \alpha_j \varphi_j$, then
\begin{align*}
H^{1/2}(\partial M) = \left\{ f \in L^2(\partial M) : \sum_{j=0}^\infty	(1 + \alpha_j)^{1/2} \left\vert\langle f, \varphi_j \rangle_{L^2(\partial M)} \right\vert^2 < \infty \right\}.
\end{align*}
We regard $H^{-1/2}(\partial M)$ as the dual space of $H^{1/2}(\partial M)$ and denote the corresponding duality pairing by $\langle \psi, f \rangle$, for $\psi \in H^{-1/2}(\partial M)$ and $f \in H^{1/2}(\partial M)$. Whenever $\psi \in L^2(\partial M)$, this pairing agrees with the usual $L^2$ pairing.

Recall that the trace operator $\operatorname{Tr} : H^1(M) \to H^{1/2}(\partial M)$ is bounded and surjective, and that $H_0^1(M) = \ker(\operatorname{Tr})$. A function $u \in H^1(M)$ is said to be \textbf{weakly harmonic} if 
\begin{align*}
\int_M\langle \nabla u, \nabla \varphi \rangle \,\dV_g = 0
\end{align*}
for every $\varphi\in H_0^1(M)$.

Given $f \in H^{1/2}(\partial M)$, there exists a unique weakly harmonic function $\hat f \in H^1(M)$ whose trace on $\partial M$ is $f$. We call $\hat f$ the harmonic extension of $f$. The Dirichlet-to-Neumann operator associated with $(M,g)$ is the bounded linear map $\mathcal D_g : H^{1/2}(\partial M)\to H^{-1/2}(\partial M)$ defined by
\begin{align*}
\langle\mathcal D_g f, \psi \rangle = \int_M \langle \nabla \hat f, \nabla \tilde \psi \rangle \, \dV_g,
\end{align*}
where $\psi\in H^{1/2}(\partial M)$ and $\tilde \psi \in H^1(M)$ is any function whose trace on $\partial M$ is $\psi$. This definition is independent of the choice of $\tilde \psi$, since the difference of two such extensions belongs to $H_0^1(M)$ and $\hat f$ is weakly harmonic.

If $f$ is smooth, elliptic regularity implies that $\hat f$ is smooth up to the boundary, and Green's formula implies that
\begin{align*}
\mathcal D_g f = \frac{\partial \hat f}{\partial \nu} = \partial_\nu \hat f
\end{align*}
in the classical sense, where $\nu$ denotes the outward unit normal to $\partial M$. This explains the terminology Dirichlet-to-Neumann operator.

The operator $\mathcal D_g$ is an elliptic pseudodifferential operator of order one. Its standard realisation as an unbounded operator on $L^2(\partial M)$, with domain $H^1(\partial M)$, is nonnegative and self-adjoint with compact resolvent. Its eigenvalues, called the Steklov eigenvalues of $(M, g)$, form a nondecreasing sequence
\begin{align*}
0 = \sigma_0(g) < \sigma_1(g) \leq \sigma_2(g) \leq \cdots \leq \sigma_k(g) \leq \cdots \nearrow + \infty,
\end{align*}
where each eigenvalue is repeated according to its multiplicity. It is well known that these eigenvalues admit the variational characterisation
\begin{align*}
\sigma_k(g) = \min_{V\in\mathcal V_{k+1}(M)} \max_{0\neq u\in V} \frac{\displaystyle \int_M \vert \nabla u \vert^2 \, \dV_g}{\displaystyle \int_{\partial M}u^2 \, \dA_g},
\end{align*}
where $\mathcal V_k(M)$ denotes the set of all $k$-dimensional subspaces of $H^1(M)$. Here and below, a Rayleigh quotient with vanishing denominator is understood to be $+\infty$.

Equivalently, a number $\sigma \geq 0$ is a Steklov eigenvalue of $(M, g)$ if there exists a nonzero function $u\in H^1(M)$ such that
\begin{align*}
\int_M \langle \nabla u, \nabla \varphi \rangle \, \dV_g = \sigma \int_{\partial M} u \varphi \, \dA_g
\end{align*}
for every $\varphi\in H^1(M)$. Standard elliptic regularity implies that such a function is smooth, and the weak formulation is equivalent to the boundary value problem
\begin{align*}
\begin{cases}
\Delta_g u = 0 & \text{in } M,\\
\partial_\nu u = \sigma u & \text{on }\partial M.
\end{cases}
\end{align*}
In particular, $u \vert_{\partial M}$ is an eigenfunction of $\mathcal D_g$ with eigenvalue $\sigma$.

We now introduce the corresponding $G$-invariant problem. Suppose that $G$ is a compact Lie group acting smoothly by isometries on $(M, g)$. If $E$ is a (Banach) space of functions on $M$ (or on $\partial M$), we consider the induced left action of $G$ on $E$ given by
\begin{align*}
(\gamma\cdot u)(p) = u(\gamma^{-1} \cdot p), \qquad \gamma\in G,
\end{align*}
for any $p \in M$ (or $p \in \partial M$). The induced action on the dual space $E^\ast$ is defined by
\begin{align*}
\langle \gamma \cdot \psi, u \rangle = \langle \psi, \gamma^{-1} \cdot u \rangle, \qquad \gamma\in G,
\end{align*}
for $\psi\in E^\ast$ and $u\in E$. We denote by
\begin{align*}
E^G = \{u \in E :\gamma \cdot u = u \text{ for every } \gamma \in G\}
\end{align*}
the corresponding (closed, if $E$ is Banach) subspace of $G$-invariant functions.

The first observation is that harmonic extension and the Dirichlet-to-Neumann operator are compatible with the $G$-action.

\begin{lem}
The Dirichlet-to-Neumann operator is $G$-equivariant, that is,
\begin{align*}
\mathcal D_g(\gamma \cdot f) = \gamma \cdot \mathcal D_g f
\end{align*}
for every $\gamma \in G$ and $f \in H^{1/2}(\partial M)$. In particular, it restricts to a bounded linear map
\begin{align*}
\mathcal D_g^G : H^{1/2}(\partial M)^G \longrightarrow H^{-1/2}(\partial M)^G.
\end{align*}
\end{lem}

\begin{proof}
Let $f\in H^{1/2}(\partial M)$ and $\gamma \in G$. We first claim that the harmonic extension of $\gamma \cdot f$ is $\gamma \cdot \hat f$. Indeed, for every $\varphi\in H_0^1(M)$, invariance of the metric and a change of variables give
\begin{align*}
\int_M \langle \nabla(\gamma \cdot \hat f), \nabla \varphi \rangle \, \dV_g = \int_M \langle \nabla \hat f, \nabla(\gamma^{-1} \cdot \varphi) \rangle \, \dV_g	= 0.
\end{align*}
Moreover, the trace of $\gamma \cdot \hat f$ is $\gamma \cdot f$. Hence, by uniqueness of harmonic extension, $\widehat{\gamma \cdot f} = \gamma \cdot \hat f$, as claimed.
	
Now let $\psi \in H^{1/2}(\partial M)$ and choose an extension $\tilde \psi \in H^1(M)$ of $\psi$. Then $\gamma^{-1} \cdot \tilde \psi$ is an extension of $\gamma^{-1} \cdot \psi$, and therefore
\begin{align*}
\langle \mathcal D_g(\gamma \cdot f), \psi \rangle &= \int_M \langle \nabla(\gamma \cdot \hat f), \nabla \tilde \psi \rangle \, \dV_g \\
&= \int_M \langle \nabla \hat f, \nabla(\gamma^{-1} \cdot \tilde \psi) \rangle \, \dV_g \\
&= \langle \mathcal D_g f, \gamma^{-1} \cdot \psi \rangle \\
&= \langle \gamma \cdot \mathcal D_g f, \psi \rangle.
\end{align*}
This concludes the proof.
\end{proof}

The $L^2$-realisation of $\mathcal D_g$ therefore restricts to a nonnegative self-adjoint operator on $L^2(\partial M)^G$ with domain $H^1(\partial M)^G$. Moreover, since $H^{1/2}(\partial M)^G$ is a closed subspace of $H^{1/2}(\partial M)$, the embedding $H^{1/2}(\partial M)^G\hookrightarrow L^2(\partial M)^G$ is compact. Consequently, the standard variational construction applies in the $G$-invariant setting. We call the eigenvalues of this restriction the \textbf{$G$-invariant Steklov eigenvalues}. Whenever the $k$-th $G$-invariant Steklov eigenvalue exists, we denote it by $\sigma_k^G(g)$. It admits the variational characterisation
\begin{align*}
\sigma_k^G(g) = \min_{V \in \mathcal V_{k+1}^G(M)} \max_{0 \neq u \in V} \frac{\displaystyle \int_M \vert \nabla u \vert^2 \, \dV_g}{\displaystyle \int_{\partial M} u^2 \, \dA_g},
\end{align*}
where $\mathcal V_k^G(M)$ denotes the set of all $k$-dimensional subspaces of $H^1(M)^G$.

Equivalently, $\sigma \geq 0$ is a $G$-invariant Steklov eigenvalue if there exists a nonzero function $u \in H^1(M)^G$ such that
\begin{align*}
\int_M \langle \nabla u, \nabla \varphi \rangle \, \dV_g = \sigma\int_{\partial M} u \varphi \, \dA_g
\end{align*}
for every $\varphi\in H^1(M)^G$. Averaging an arbitrary test function over $G$ shows that the same identity then holds for every $\varphi\in H^1(M)$. Thus a $G$-invariant Steklov eigenfunction is simply a $G$-invariant solution of the usual Steklov problem.

In general, the $G$-invariant Steklov spectrum need not be infinite, as mentioned in the Introduction. The following lemma characterises when the space of smooth $G$-invariant functions on the boundary is infinite-dimensional.

\begin{lem}\label{invariant functions}
Let $N$ be a compact smooth manifold, possibly disconnected, and let $G$ be a compact Lie group acting smoothly on $N$. For each connected component $\Sigma$ of $N$, let
\begin{align*}
G_\Sigma = \{\gamma \in G : \gamma(\Sigma) = \Sigma\}.
\end{align*}
Then $C^\infty(N)^G$ is infinite-dimensional if and only if there exists a connected component $\Sigma$ of $N$ such that the action of $G_\Sigma$ on $\Sigma$ is nontransitive.
\end{lem}

\begin{proof}
Since $N$ has finitely many connected components, $G$ permutes them. Let $\Sigma_1, \dots, \Sigma_r$ be representatives of the $G$-orbits of connected components, and set $N_i = G(\Sigma_i)$. Then
\begin{align*}
N = N_1\sqcup\cdots\sqcup N_r.
\end{align*}
The restriction map
\begin{align*}
C^\infty(N_i)^G & \longrightarrow C^\infty(\Sigma_i)^{G_{\Sigma_i}} \\
f & \longmapsto f \vert_{\Sigma_i},
\end{align*}
is an isomorphism of vector spaces. Thus
\begin{align*}
C^\infty(N)^G \cong \bigoplus_{i = 1}^r C^\infty(\Sigma_i)^{G_{\Sigma_i}}.
\end{align*}
	
It is therefore enough to consider the case in which $N$ is connected. If the action is transitive, every $G$-invariant smooth function is constant.
	
Assume now that the action is nontransitive. Choose $p, q \in N$ belonging to distinct $G$-orbits. Since $G$ is compact, the orbits $G \cdot p$ and $G \cdot q$ are disjoint compact subsets of $N$. Hence there exists $\xi \in C^\infty(N)$ such that
\begin{align*}
\xi \equiv 0 \text{ on } G \cdot p, \qquad \xi \equiv 1 \text{ on } G\cdot q.
\end{align*}
Let $\dd \mu$ be the normalised Haar measure on $G$, and define
\begin{align*}
f(x) = \int_G \xi(\gamma \cdot x) \, \dmu(\gamma).
\end{align*}
Then $f$ is smooth and $G$-invariant. Moreover, $f(p) = 0$ and $f(q) = 1$, so $f$ is nonconstant. Since $N$ is connected, $f(N)$ is an interval containing more than one point. If $a_0, \dots, a_m \in \mathbb R$ satisfy
\begin{align*}
a_0 + a_1 f + \cdots + a_m f^m = 0,
\end{align*}
then the polynomial $P(t) = a_0 + a_1 t + \cdots + a_m t^m$ vanishes on $f(N)$ and hence is identically zero. Thus $1, f, f^2, \dots$ are linearly independent, and $C^\infty(N)^G$ is infinite-dimensional.
\end{proof}

Combining Lemma \ref{invariant functions} with the spectral properties discussed above, we obtain the following.

\begin{prop}
Let $(M, g)$ be a compact connected smooth Riemannian manifold with $\partial M \neq \emptyset$, and let $G$ be a compact Lie group acting smoothly by isometries on $M$. Assume that the action is boundary-nontransitive. Then the $G$-invariant Steklov spectrum is an infinite discrete sequence of eigenvalues of finite multiplicity,
\begin{align*}
0 = \sigma_0^G(g) < \sigma_1^G(g) \leq \sigma_2^G(g) \leq \cdots \leq \sigma_k^G(g) \leq \cdots \nearrow + \infty.
\end{align*}
\end{prop}

\begin{proof}	
By Lemma \ref{invariant functions}, the space $C^\infty(\partial M)^G$, and hence $L^2(\partial M)^G$, is infinite-dimensional. Since the $L^2$-realisation of $\mathcal D_g^G$ is nonnegative and self-adjoint with compact resolvent, its spectrum consists of an infinite sequence of eigenvalues of finite multiplicity diverging to $+\infty$. Finally, since $M$ is connected, the kernel consists of the constant functions, so $0$ is simple.
\end{proof}

\begin{remark} \label{remark: n geq 3 proof}
We justify the assertion in Remark \ref{remark: n geq 3 in main theorem}. Let $G$ be a compact Lie group with $\dim G \geq 1$ acting smoothly and effectively on a compact connected smooth surface $M$ with nonempty boundary. Suppose, for a contradiction, that there exists a boundary component $N \subseteq \partial M$ on which
$G_N$ acts nontransitively. Since the identity component $G^0$ is connected, it preserves every connected component of $\partial M$, so $G^0 \subseteq G_N$. In particular, $G^0$ cannot act transitively on $N$. Since $N$ is a circle, any nontrivial action of the compact and connected group $G^0$ on $N$ is transitive. Hence $G^0$ fixes $N$ pointwise.

Now choose any $G$-invariant Riemannian metric on $M$. Fix $p \in N$ and $\gamma \in G^0$. Since $\gamma$ fixes $N$ pointwise, $\dd \gamma_p \vert_{T_p N} = \id_{T_p N}$. Moreover, since $\gamma$ is an isometry of $M$, $\dd \gamma_p$ preserves the inward unit normal. Thus $\dd \gamma_p = \id_{T_p M}$. An isometry of a connected Riemannian manifold is determined by its value and differential at one point, so $\gamma = \id_M$. Since $\gamma\in G^0$ was arbitrary, $G^0$ acts trivially on $M$. Effectiveness then implies $G^0 = \{ e \}$, contradicting $\dim G^0 = \dim G \geq 1$.
\end{remark}

\subsection{Steklov--Neumann problem} \label{sec: Steklov-Neumann problem}

We now recall the Steklov--Neumann problem, also known as the sloshing problem, and introduce its $G$-invariant counterpart. This problem will play an important role in the proof of Theorem \ref{thm: infinite sup}.

Let $(\Omega,g)$ be a compact connected Riemannian manifold of dimension $n\geq2$ with Lipschitz boundary. Suppose that
\begin{align*}
\partial \Omega = \overline{\partial_V \Omega} \cup \overline{\partial_L \Omega},
\end{align*}
where $\partial_V \Omega$ and $\partial_L \Omega$ are disjoint relatively open smooth subsets of $\partial\Omega$, with $\partial_V\Omega\neq\emptyset$, and whose closures meet along a smooth codimension-two submanifold. On $\partial_V \Omega$ we impose the Steklov boundary condition, while on $\partial_L\Omega$ we impose the Neumann boundary condition.

A number $\mu \geq 0$ is called a Steklov--Neumann eigenvalue of $(\Omega, g)$ if there exists a nonzero function $u \in H^1(\Omega)$ such that
\begin{align*}
\int_\Omega \langle \nabla u, \nabla \varphi \rangle \, \dV_g =	\mu \int_{\partial_V \Omega} u \varphi \, \dA_g
\end{align*}
for every $\varphi \in H^1(\Omega)$. Such a function is called a Steklov--Neumann eigenfunction. Taking test functions in $H_0^1(\Omega)$ shows that $u$ is weakly harmonic in $\Omega$. Moreover, elliptic regularity implies that $u$ is smooth in $\Omega$ and up to the relative interior of each of the boundary pieces $\partial_V \Omega$ and $\partial_L \Omega$. There it satisfies
\begin{align*}
\begin{cases}
\Delta_g u = 0 & \text{in } \Omega,\\
\partial_\nu u = \mu u & \text{on } \partial_V \Omega,\\
\partial_\nu u = 0 & \text{on } \partial_L \Omega.
\end{cases}
\end{align*}

By the compactness of the trace embedding $H^1(\Omega)\hookrightarrow L^2(\partial_V\Omega)$, the Steklov--Neumann eigenvalues form a nondecreasing sequence
\begin{align*}
0 = \mu_0(\Omega, g) < \mu_1(\Omega, g) \leq \mu_2(\Omega, g) \leq \cdots \leq \mu_k(\Omega, g) \leq \cdots \nearrow + \infty,
\end{align*}
where each eigenvalue is repeated according to its multiplicity. Since $\Omega$ is connected, the eigenspace corresponding to $\mu_0(\Omega, g) = 0$ consists precisely of the constant functions. It is well known that these eigenvalues admit the variational characterisation
\begin{align*}
\mu_k(\Omega, g) = \min_{V \in \mathcal V_{k+1}(\Omega)} \max_{0 \neq u \in V} 	\frac{\displaystyle \int_\Omega \vert \nabla u \vert^2 \, \dV_g}{\displaystyle \int_{\partial_V \Omega} u^2 \, \dA_g},
\end{align*}
where $\mathcal V_k(\Omega)$ denotes the set of all $k$-dimensional subspaces of $H^1(\Omega)$. 

Suppose now that $G$ is a compact Lie group acting smoothly by isometries on $(\Omega,g)$ and preserving both $\partial_V\Omega$ and $\partial_L\Omega$. Restricting the preceding variational problem to $G$-invariant functions gives the $G$-invariant Steklov--Neumann problem. Whenever the $k$-th $G$-invariant Steklov--Neumann eigenvalue exists, we denote it by $\mu_k^G(\Omega,g)$. It admits the variational characterisation
\begin{align*}
\mu_k^G(\Omega, g) = \min_{V\in\mathcal V_{k+1}^G(\Omega)} \max_{0 \neq u\in V}	\frac{\displaystyle \int_\Omega \vert \nabla u \vert^2 \, \dV_g}{\displaystyle \int_{\partial_V \Omega} u^2 \, \dA_g},
\end{align*}
where $\mathcal V_k^G(\Omega)$ denotes the set of all $k$-dimensional subspaces of $H^1(\Omega)^G$. As in the Steklov case, every $G$-invariant Steklov--Neumann eigenfunction is a $G$-invariant weak solution of the usual Steklov--Neumann problem.

In general, the $G$-invariant Steklov--Neumann spectrum need not be infinite. In what follows, we only use $\mu_k^G(\Omega,g)$ when the corresponding eigenvalue exists.

\begin{remark} \label{rem: disconnected problem}
The spectral constructions and variational characterisations in this section also apply without assuming that $M$ or $\Omega$ is connected. For the Steklov problem, it is enough to assume that every connected component of $M$ has nonempty boundary; for the Steklov--Neumann problem, it is enough to assume that every connected component of $\Omega$ meets $\partial_V \Omega$.

Under these assumptions, the problem on the whole manifold decomposes into the corresponding problems on its finitely many connected components. Connectedness is therefore needed only to ensure the simplicity of the zero eigenvalue in the ordinary spectrum. Indeed, the zero eigenspace consists precisely of the functions that are constant on each connected component, so its dimension equals the number of components. For the $G$-invariant problem, these constants must agree on components belonging to the same orbit of the induced action of $G$ on the set of connected components. Thus the multiplicity of zero in the invariant spectrum equals the number of such orbits. In particular, zero remains simple whenever $G$ permutes the connected components transitively, even if the manifold itself is disconnected. As before, we list eigenvalues with multiplicity, including zero.
\end{remark}

\section{Extremal metrics for the Steklov eigenvalues} \label{sec: Extremal metrics for the Steklov eigenvalues}

In this section, we introduce the notion of extremal metrics for the ($G$-invariant) Steklov eigenvalues, inspired by the works of Nadirashvili \cite{nadirashvili}, Karpukhin--M\'etras \cite{karpukhin2022} and Colbois--Dryden--El Soufi \cite{colbois_G}. We then follow the techniques of Fraser--Schoen \cite{fraser} and Karpukhin--M\'etras to compute the derivatives of these eigenvalues without invoking Kato's analytic perturbation theory, while also correcting some imprecisions in the literature.

Let $M$ be a compact connected smooth manifold of dimension $n \geq 2$ with nonempty boundary. For each smooth Riemannian metric $g$ on $M$, denote by $\sigma_k(g)$ the $k$-th Steklov eigenvalue and by $E_k(g) \subset L^2(\partial M, g)$ the corresponding eigenspace. We consider the normalised eigenvalue
\begin{align*}
\overline{\sigma}_k(g) = N(g) \sigma_k(g), \qquad N(g) = \Vol(M, g)^{\frac{2 - n}{n}} \Vol(\partial M, g).
\end{align*}
The reader may naturally wonder about the origin of this normalisation. Karpukhin--M\'etras showed in \cite[Proposition 1]{karpukhin2022} that no other scale-invariant normalisation involving powers of the interior and boundary volumes admits extremal metrics.

For $k \geq 1$, we say that a smooth Riemannian metric $g$ is \textbf{$\overline{\sigma}_k$-extremal} if, for every smooth $1$-parameter family of Riemannian metrics $g(t)$ defined near $t = 0$ and satisfying $g(0) = g$, either
\begin{align*}
\overline{\sigma}_k(g(t)) \leq \overline{\sigma}_k(g) + o(t) \quad \text{or} \quad \overline{\sigma}_k(g(t)) \geq \overline{\sigma}_k(g) + o(t)
\end{align*}
as $t \to 0$. We say that $g$ is \textbf{$\overline{\sigma}_k$-conformally extremal} if this condition holds for every such family contained in the conformal class $[g]$.

Now fix a smooth effective action of a compact Lie group $G$ on $M$, and assume that the action is boundary-nontransitive. For a smooth $G$-invariant Riemannian metric $g$, let $\sigma_k^G(g)$ denote the $k$-th $G$-invariant Steklov eigenvalue and set
\begin{align*}
\overline{\sigma}_k^G(g) = N(g) \sigma_k^G(g).
\end{align*}
For $k \geq 1$, we say that $g$ is \textbf{$\overline{\sigma}_k^G$-extremal} if, for every smooth $1$-parameter family of $G$-invariant Riemannian metrics $g(t)$ defined near $t = 0$ and satisfying $g(0) = g$, either
\begin{align*}
\overline{\sigma}_k^G(g(t)) \leq \overline{\sigma}_k^G(g) + o(t) \quad \text{or} \quad \overline{\sigma}_k^G(g(t)) \geq \overline{\sigma}_k^G(g) + o(t)
\end{align*}
as $t \to 0$. We say that $g$ is \textbf{$\overline{\sigma}_k^G$-conformally extremal} if this condition holds for every such family contained in $[g]$.

\subsection{First variation of the Steklov eigenvalues} \label{sec: First variation of the Steklov eigenvalues}

The Steklov eigenvalue functional $\sigma_k$ is locally Lipschitz on the space of Riemannian metrics endowed with the $C^\infty$ topology; see, for instance, \cite{fraser,kriegl}. The next proposition gives a formula for its derivative along a smooth $1$-parameter family of metrics.

\begin{prop} \label{derivative sigma}
Let $(-\varepsilon, \varepsilon) \ni t \mapsto g(t)$ be a smooth $1$-parameter family of Riemannian metrics on $M$ and write $\sigma_k(t) := \sigma_k(g(t))$. If $\sigma_{k-1}(t_0) < \sigma_k(t_0)$ or $\sigma_{k}(t_0) < \sigma_{k+1}(t_0)$, and $\dot{\sigma_k}(t_0)$ exists, then
\begin{align*}
\dot{\sigma}_k(t_0) = - \int_M \langle q_{g(t_0)}(\hat u), h \rangle \dV_{g(t_0)} - \frac{\sigma_k(t_0)}{2} \int_{\partial M} u^2 \langle g(t_0)\vert_{\partial M}, h\vert_{\partial M} \rangle \dA_{g(t_0)},
\end{align*}
for any $u \in E_k(g(t_0))$ with $\norm{u}_{L^2(\partial M, g(t_0))} = 1$. Here, $h = \left. \ddt \right\vert_{t = t_0} g(t)$, $\langle \cdot, \cdot \rangle$ is the inner product induced by $g(t_0)$ in the space $S^2(M)$ (or $S^2(\partial M)$) of smooth covariant $2$-tensors, $\dV_{g(t_0)}$ and $\dA_{g(t_0)}$ are the volume elements of $M$ and $\partial M$ induced by $g(t_0)$, and $q_g : C^{\infty}(M) \to S^2(M)$ is given by
\begin{align*}
q_g(f) = \dd f \otimes \dd f - \frac{1}{2} \vert \nabla f \vert^2 g.
\end{align*}
\end{prop}

\begin{proof}
Suppose $\dot{\sigma}_k(t_0)$ exists and $\sigma_{k-1}(t_0) < \sigma_k(t_0)$; the other inequality is treated similarly. Choose a neighbourhood of $t_0$ where this inequality remains valid and define on this neighbourhood
\begin{align*}
\mathcal{E}(t) = \sum_{i=0}^{k-1} E_i(g(t)) \subset L^2(\partial M, g(t)).
\end{align*} 
Then $\mathcal{E}(t)$ has constant dimension $k$. Moreover, $\mathcal{E}(t)$ varies smoothly in $t$ for all $t$ close to $t_0$. Let $P_t : L^2(\partial M, g(t)) \to \mathcal{E}(t)$ be the orthogonal projection onto $\mathcal{E}(t)$. Choose $u \in E_k(g(t_0))$ with $\norm{u}_{L^2(\partial M, g(t_0))} = 1$ and let $u_t = u - P_t(u)$. Since $P_{t_0}(u) = 0$, we have $u_{t_0} = u$. Let $\hat u_t$ denote the harmonic extension of $u_t$ with respect to $g(t)$, and write
\begin{align*}
\dot u_{t_0} = \left. \ddt \right\vert_{t = t_0} u_t, \qquad \dot{\hat u}_{t_0} = \left. \ddt \right\vert_{t = t_0} \hat u_t.
\end{align*}
Then the real function
\begin{align*}
f(t) = \int_M g(t)(\nabla \hat u_t, \nabla \hat u_t) \dV_{g(t)} - \sigma_k(t) \int_{\partial M} u_t^2 \dA_{g(t)}
\end{align*}
satisfies $f(t) \geq 0$ in a neighbourhood of $t_0$ and $f(t_0) = 0$, so that $\dot f(t_0) = 0$. It remains to compute this derivative. Recall that
\begin{align*}
\left. \ddt \right\vert_{t = t_0} \dV_{g(t)} &= \frac{1}{2} \langle g(t_0), h \rangle \dV_{g(t_0)}, \\
\left. \ddt \right\vert_{t = t_0} \dA_{g(t)} &= \frac{1}{2} \langle g(t_0)\vert_{\partial M}, h\vert_{\partial M} \rangle \dA_{g(t_0)}.
\end{align*}
Moreover,
\begin{align*}
\left. \ddt \right\vert_{t = t_0} g(t)(\nabla \hat u_t, \nabla \hat u_t) = 2 g(t_0)(\nabla \hat u, \nabla \dot{\hat u}_{t_0}) - h(\nabla \hat u, \nabla \hat u).
\end{align*}
Consequently,
\begin{align*}
\dot f(t_0) &= \int_M \left( \left. \ddt \right\vert_{t = t_0} g(t)(\nabla \hat u_t, \nabla \hat u_t) \right) \dV_{g(t_0)} + \int_M g(t_0)(\nabla \hat u, \nabla \hat u) \left( \left. \ddt \right\vert_{t = t_0} \dV_{g(t)} \right) \\
&\quad - \dot{\sigma}_k(t_0) - \sigma_k(t_0) \int_{\partial M} \left( \left. \ddt \right\vert_{t = t_0} u_t^2 \right) \dA_{g(t_0)} - \sigma_k(t_0) \int_{\partial M} u^2 \left( \left. \ddt \right\vert_{t = t_0} \dA_{g(t)} \right) \\
&= \int_M \left[ 2 g(t_0)(\nabla \hat u, \nabla \dot{\hat u}_{t_0}) - h(\nabla \hat u, \nabla \hat u) + \frac{1}{2} \vert \nabla \hat u \vert^2 \langle g(t_0), h \rangle \right] \dV_{g(t_0)} \\
&\quad - \dot{\sigma}_k(t_0) - \sigma_k(t_0) \int_{\partial M} \left( 2u \dot u_{t_0} + \frac{1}{2} u^2 \langle g(t_0)\vert_{\partial M}, h\vert_{\partial M} \rangle \right) \dA_{g(t_0)} \\
&= \int_M \left[ 2 g(t_0)(\nabla \hat u, \nabla \dot{\hat u}_{t_0}) - \left\langle \dd \hat u \otimes \dd \hat u - \frac{1}{2} \vert \nabla \hat u \vert^2 g(t_0), h \right\rangle \right] \dV_{g(t_0)} \\
&\quad - \dot{\sigma}_k(t_0) - \sigma_k(t_0) \int_{\partial M} \left( 2u \dot u_{t_0} + \frac{1}{2} u^2 \langle g(t_0)\vert_{\partial M}, h\vert_{\partial M} \rangle \right) \dA_{g(t_0)}.
\end{align*}
Since $u \in E_k(g(t_0))$ and the boundary trace of $\dot{\hat u}_{t_0}$ is $\dot u_{t_0}$, we also have
\begin{align*}
\int_M g(t_0)(\nabla \hat u, \nabla \dot{\hat u}_{t_0}) \dV_{g(t_0)} = \sigma_k(t_0) \int_{\partial M} u \dot u_{t_0} \dA_{g(t_0)}.
\end{align*}
Thus, using $\dot f(t_0) = 0$ and cancelling the terms involving $\dot u_{t_0}$, we obtain
\begin{align*}
\dot{\sigma}_k(t_0) = - \int_M \left\langle \dd \hat u \otimes \dd \hat u - \frac{1}{2} \vert \nabla \hat u \vert^2 g(t_0), h \right\rangle \dV_{g(t_0)} - \frac{\sigma_k(t_0)}{2} \int_{\partial M} u^2 \langle g(t_0)\vert_{\partial M}, h\vert_{\partial M} \rangle \dA_{g(t_0)},
\end{align*}
as desired.

\end{proof}

We now specialize the preceding analysis to deformations within a fixed conformal class. Let $g(t) =	e^{\varphi(t)}g$ be a smooth conformal deformation of a fixed Riemannian metric $g$. At any $t_0$, its infinitesimal variation is $h = \ddt g(t) \vert_{t = t_0} = \dot\varphi(t_0) g(t_0)$. Consequently,
\begin{align*}
\langle q_{g(t_0)}(\hat u), h \rangle = \frac{(2-n)}{2}	\vert\nabla\hat u\vert^2 \dot{\varphi}(t_0)
\end{align*}
on $M$, while
\begin{align*}
\langle g(t_0) \vert_{\partial M}, h \vert_{\partial M}	\rangle = (n - 1) \dot{\varphi}(t_0)
\end{align*}
on $\partial M$. Substituting these identities into the general first-variation formula of Proposition \ref{derivative sigma} yields the following conformal version.

\begin{prop}\label{derivative sigma conformal}
Let $(-\varepsilon, \varepsilon) \ni t \mapsto g(t) = e^{\varphi(t)} g$	be a smooth $1$-parameter family of Riemannian metrics on $M$ conformal to a fixed metric $g$. Suppose that either $\sigma_{k-1}(t_0) < \sigma_k(t_0)$ or $\sigma_k(t_0) < \sigma_{k + 1}(t_0)$, and that $\dot{\sigma}_k(t_0)$ exists. Then
\begin{align*}
\dot{\sigma}_k(t_0) = \frac{(n - 2)}{2} \int_M \vert \nabla \hat u \vert^2 \dot\varphi(t_0) \dV_{g(t_0)} - \frac{(n - 1)}{2} \sigma_k(t_0) \int_{\partial M} u^2 \dot \varphi(t_0) \dA_{g(t_0)}
\end{align*}
for every $u \in E_k(g(t_0))$ satisfying $\norm{u}_{L^2(\partial M, g(t_0))} = 1$.
\end{prop}

\subsection{Extremality under general metric deformations} \label{sec: Extremality under general metric deformations}

We now use the first-variation formula established in the preceding subsection to derive a consequence of the extremality of a metric. The main point is that, after fixing the normalisation, the admissible infinitesimal deformations form a codimension-one subspace, and extremality forces the quadratic form associated with the variation of the eigenvalue to vanish on some normalised eigenfunction in every such direction.

We first describe the infinitesimal constraint imposed by fixing the normalisation and introduce the corresponding quadratic forms for general interior and boundary tensor directions.

Given a Riemannian metric $g$ on $M$, let $L^2(S^2(M),g)$ denote the Hilbert space of $L^2$ symmetric covariant $2$-tensors on $M$, endowed with the inner product
\begin{align*}
\langle\!\langle H_1, H_2 \rangle\!\rangle = \int_M \langle H_1, H_2 \rangle \dV_g.
\end{align*}
We define $L^2(S^2(\partial M),g\vert_{\partial M})$ analogously. Let $\mathcal{H} = L^2(S^2(M),g) \times L^2(S^2(\partial M), g\vert_{\partial M})$ be the corresponding product Hilbert space, endowed with the inner product
\begin{align*}
\langle\!\langle (H_1, \eta_1), (H_2, \eta_2) \rangle\!\rangle_{\mathcal{H}} = \langle\!\langle H_1, H_2 \rangle\!\rangle + \langle\!\langle \eta_1, \eta_2 \rangle\!\rangle.
\end{align*}

For a pair $(H, \eta) \in \mathcal{H}$ and an integer $s \geq 1$, define a quadratic form $Q_{(H,\eta)}^{g,s} : C^\infty(\partial M) \to \mathbb{R}$ by
\begin{align*}
Q_{(H, \eta)}^{g,s}(u) = - \int_M \langle q_g(\hat u), H \rangle \dV_{g} - \frac{\sigma_s(g)}{2} \int_{\partial M} u^2 \langle g\vert_{\partial M}, \eta \rangle \dA_{g},
\end{align*}
where $q_g$ is defined in Proposition \ref{derivative sigma}. With this notation, the derivative formula in that proposition can be written as
\begin{align*}
\dot \sigma_k(t_0) = Q_{(h, h\vert_{\partial M})}^{g(t_0),k}(u), 
\end{align*}
for any $u \in E_k(g(t_0))$ satisfying $\norm{u}_{L^2(\partial M,g(t_0))} = 1$.

Since the normalised eigenvalue satisfies $\overline{\sigma}_k(g) =	N(g) \sigma_k(g)$, we will also need the first variation of the normalisation factor $N$. The following lemma gives this variation and identifies its gradient with respect to the product Hilbert space $\mathcal H$ introduced above.

\begin{lem} \label{normalisation}
Given $(-\varepsilon, \varepsilon) \ni t \mapsto g(t)$ a smooth $1$-parameter family of Riemannian metrics on $M$ with $g(0) = g$ and $\ddt g(t) \vert_{t=0} = h$, we have:
\begin{align*}
\left. \ddt \right\vert_{t = 0} N(g(t)) = \frac{N(g)}{2} \left[ \frac{2 - n}{n \Vol(M,g)} \int_M \langle g, h \rangle \dV_g + \frac{1}{\Vol(\partial M,g)} \int_{\partial M} \langle g\vert_{\partial M}, h\vert_{\partial M} \rangle \dA_g \right].
\end{align*}
Alternatively, regarding $N$ as a function of the pair $(g, g\vert_{\partial M})$, we can write
\begin{align*}
\nabla N(g, g\vert_{\partial M}) = \frac{N(g)}{2} \left( \frac{2 - n}{n \Vol(M,g)} \, g, \frac{1}{\Vol(\partial M,g)} \, g\vert_{\partial M} \right),
\end{align*}
so that
\begin{align*}
\left.\ddt\right|_{t=0} N(g(t)) = \langle\!\langle \nabla N(g, g \vert_{\partial M}), (h, h \vert_{\partial M}) \rangle\!\rangle_{\mathcal{H}}.
\end{align*}
\end{lem}

The preceding lemma identifies the tangent space at $(g,g\vert_{\partial M})$ to the level set of the normalisation functional $N$. More precisely, a pair $(H,\eta)\in\mathcal H$ is tangent to this level set precisely when
\begin{align*}
\left\langle\!\left\langle	\nabla N(g, g\vert_{\partial M}), (H,\eta) \right\rangle\!\right\rangle_{\mathcal H} =	0.
\end{align*}
On the other hand, infinitesimal variations of smooth Riemannian metrics give rise only to pairs of the form $ (h,h|_{\partial M})$, for $h\in S^2(M)$. To pass from such genuine metric variations to arbitrary elements of $\mathcal H$, we will use the following approximation result.

\begin{lem} \label{denseness}
The subspace
\begin{align*}
S = \{ (h, h\vert_{\partial M}) \in \mathcal{H} : h \in S^2(M) \}
\end{align*}
is dense in $\mathcal{H}$.
\end{lem}

\begin{proof}

Let $(H, \eta) \in \mathcal{H}$ be arbitrary. Since $S^2(M)$ is dense in $L^2(S^2(M), g)$ and $S^2(\partial M)$ is dense in $L^2(S^2(\partial M), g\vert_{\partial M})$, there are sequences $(H_i) \subset S^2(M)$ and $(\eta_i) \subset S^2(\partial M)$ such that $H_i \to H$ and $\eta_i \to \eta$ in the corresponding $L^2$ norms, as $i \to \infty$. We cannot conclude the result yet, as the pair $(H_i, \eta_i)$ does not belong to $S$. The idea is to modify $H_i$ in a thin collar around $\partial M$ so that its restriction to $\partial M$ coincides with $\eta_i$; making the thickness of this collar converge to zero as $i \to \infty$ guarantees that the modified $H_i$ still converges to $H$.

Let $\Phi : \partial M \times [0, \varepsilon] \to U \subset M$ be a smooth collar diffeomorphism onto a compact neighbourhood $U$ of $\partial M$, with $\Phi(p, 0) = p$ for all $p \in \partial M$. Let $\pi : U \to \partial M$ be the projection defined by $\pi(\Phi(p, t)) = p$. 

Now let $a_i = \eta_i - H_i\vert_{\partial M}$ and $\tilde{a}_i = \pi^\ast a_i \in S^2(U)$. Pick a cutoff function $\chi \in C^\infty(\mathbb{R})$ with $\chi(0) = 1$ and whose support is contained in $[-1,1]$. For $\delta \in (0, \varepsilon]$, set
\begin{align*}
R_{i, \delta}(\Phi(p, t)) = \chi \left( \frac{t}{\delta} \right)\tilde{a}_i(\Phi(p, t))
\end{align*}
on $U$ and $R_{i, \delta} \equiv 0$ on $M \setminus U$, so that $R_{i, \delta} \in S^2(M)$. Notice that $R_{i, \delta}$ coincides with $a_i$ along $\partial M$, for if $p \in \partial M$,
\begin{align*}
R_{i, \delta}(p) = R_{i, \delta}(\Phi(p, 0)) = \chi(0) \tilde{a}_i(\Phi(p, 0)) = a_i(p).
\end{align*}
Moreover, for each fixed $i\geq 1$, we have
\begin{align*}
\|R_{i,\delta}\|_{L^2(S^2(M),g)} \to 0 \qquad \text{ as }\delta \to 0.
\end{align*}
Indeed, $\tilde{a}_i$ is smooth on the compact collar $U$, hence bounded, and $\vert R_{i,\delta} \vert^2 \leq \norm{\chi}_{L^\infty(\mathbb{R})}^2 \vert \tilde{a}_i \vert^2$. Since $\operatorname{supp}(R_{i,\delta})\subset \Phi(\partial M\times[0,\delta])$ and the volume of this collar tends to zero as $\delta\to0$, the conclusion follows. Thus, for each $i\ge1$ we may choose $\delta_i\in(0,\varepsilon)$ such that $\delta_i < \frac{1}{i}$ and $\norm{R_{i,\delta_i}}_{L^2(S^2(M),g)} < \frac{1}{i}$.

Now set $h_i = H_i + R_{i, \delta_i} \in S^2(M)$. By construction, $h_i$ coincides with $\eta_i$ along $\partial M$, so the pair $(h_i, \eta_i)$ lies in $S$. Finally, 
\begin{align*}
\norm{h_i - H}_{L^2(S^2(M),g)} \leq \norm{H_i - H}_{L^2(S^2(M),g)} + \norm{R_{i, \delta_i}}_{L^2(S^2(M),g)} \to 0
\end{align*} 
and
\begin{align*}
\norm{h_i\vert_{\partial M} - \eta}_{L^2(S^2(\partial M),g\vert_{\partial M})} = \norm{\eta_i - \eta}_{L^2(S^2(\partial M),g\vert_{\partial M})} \to 0.
\end{align*}
Hence, $(h_i, h_i\vert_{\partial M}) \to (H, \eta)$ in $\mathcal{H}$.
\end{proof}

Combining the first-variation formula for $\sigma_k$, the characterisation of the tangent space to the level set of $N$, and the preceding density lemma, we obtain the following consequence of the extremality of $g$. It asserts that, in every direction tangent to the normalisation constraint, the quadratic form associated with the variation of $\sigma_k$ vanishes on some normalised eigenfunction.

\begin{prop} \label{isotropic}
If a metric $g$ is $\overline{\sigma}_k$-extremal, $k \geq 1$, then for any pair $(H, \eta) \in \mathcal{H}$ which satisfies
\begin{align*} 
\langle\!\langle \nabla N(g, g\vert_{\partial M}), (H, \eta) \rangle\!\rangle_{\mathcal{H}} = 0
\end{align*}
there exists $u \in E_k(g)$ with $\norm{u}_{L^2(\partial M,g)} = 1$ such that $Q_{(H, \eta)}^{g,k}(u) = 0$.
\end{prop}

\begin{proof}
Let $\varphi : \mathcal{H} \to \mathbb{R}$ be the continuous linear functional defined by 
\begin{align*}
\varphi(H, \eta) = \langle\!\langle \nabla N(g, g\vert_{\partial M}), (H, \eta) \rangle\!\rangle,
\end{align*}
and let $(H, \eta)$ be as in the statement, i.e., $(H, \eta) \in \ker \varphi$. By Lemma \ref{denseness}, there exists a sequence $(\tilde{h}_j) \subset S^2(M)$ such that $(\tilde{h}_j, \tilde{h}_j\vert_{\partial M}) \to (H, \eta)$ in $\mathcal{H}$. Now project onto $\ker \varphi$: for each $j \geq 1$, let 
\begin{align*}
h_j = \tilde{h}_j - \frac{\varphi(\tilde{h}_j, \tilde{h}_j\vert_{\partial M})}{\varphi(g, g\vert_{\partial M})} \, g. 
\end{align*}
Then $\varphi(h_j, h_j\vert_{\partial M}) = 0$ and 
\begin{align*}
\norm{(h_j, h_j\vert_{\partial M}) - (H, \eta)}_{\mathcal{H}} \leq \norm{(\tilde{h}_j, \tilde{h}_j\vert_{\partial M}) - (H, \eta)}_{\mathcal{H}} + \left\lvert \frac{\varphi(\tilde{h}_j, \tilde{h}_j\vert_{\partial M})}{\varphi(g, g\vert_{\partial M})} \right\rvert \norm{(g, g\vert_{\partial M})}_{\mathcal{H}} \to 0
\end{align*}
since $\varphi(\tilde{h}_j, \tilde{h}_j\vert_{\partial M}) \to \varphi(H, \eta) = 0$. Thus $(h_j, h_j\vert_{\partial M}) \in \ker \varphi$ and the sequence converges to $(H, \eta)$.

Now for each $j \geq 1$, let $(-a_j, a_j) \ni t \mapsto g_j(t)$ be the smooth $1$-parameter family of metrics defined by
\begin{align*}
g_j(t) = \left[ \frac{N(g)}{N(g+th_j)} \right]^2 (g + th_j),
\end{align*}
for $t \in (-a_j, a_j)$ and sufficiently small $a_j > 0$. Then $N(g_j(t)) = N(g)$ for all $t \in (-a_j, a_j)$, $g_j(0) = g$, and $\ddt g_j(t)|_{t = 0} = h_j$ for all $j \geq 1$. Since
\begin{align*}
\overline{\sigma}_k(g_j(t)) = N(g_j(t)) \sigma_k(g_j(t)) =	N(g)\,\sigma_k(g_j(t)),
\end{align*}
the extremality of $g$ implies that either
\begin{align*}
\sigma_k(g_j(t)) \leq \sigma_k(g) + o(t) \quad \text{ or } \quad \sigma_k(g_j(t)) \geq \sigma_k(g) + o(t)
\end{align*}
as $t \to 0$. We consider the former case; the latter is treated analogously. Choose an index $1 \leq l \leq k$ such that $\sigma_l(g) = \sigma_k(g)$ and $\sigma_{l-1}(g) < \sigma_l(g)$. Then, after possibly shrinking the interval $(-a_j, a_j)$, we have $\sigma_{l-1}(g_j(t)) < \sigma_l(g_j(t))$ for all $t \in (-a_j, a_j)$ and
\begin{align*}
\sigma_l(g_j(t)) \leq \sigma_k(g_j(t)) \leq \sigma_k(g) + o(t) = \sigma_l(g) + o(t)
\end{align*}
as $t \to 0$. In particular, taking the limit $t \to 0^-$ yields
\begin{align*}
\liminf_{t \to 0^-} \frac{\sigma_l(g_j(t)) - \sigma_l(g)}{t} \geq 0.
\end{align*}
So, there exists a sequence $\varepsilon_i > 0$ decreasing to $0$ such that
\begin{align*}
-\frac{1}{i} < \frac{\sigma_l(g_j(-\varepsilon_i)) - \sigma_l(g)}{-\varepsilon_i} = \frac{1}{\varepsilon_i} \int_{-\varepsilon_i}^0 \dot{\sigma}_l(g_j(t)) \, \mathrm{d}t \leq \operatorname{ess \, sup} \{   \dot{\sigma}_l(t) : t \in [-\varepsilon_i, 0] \}. 
\end{align*}
We can thus find a sequence of $t_i < 0$ increasing to $0$ such that $\dot{\sigma}_l(g_j(t_i))$ exists and is given by 
\begin{align*}
\dot{\sigma}_l(g_j(t_i)) = Q_{\omega_{ij}}^{g_j(t_i), l}(u^{(j)}_i) > -\frac{1}{i},
\end{align*}
where $\omega_{ij} = (\dot g_j(t_i), \dot g_j(t_i) \vert_{\partial M}) \in \mathcal{H}$ and $u^{(j)}_i \in E_l(g_j(t_i))$ with $\norm{u^{(j)}_i}_{L^2(\partial M, g_j(t_i))} = 1$ (see Proposition \ref{derivative sigma}). Passing to a subsequence if necessary, $u^{(j)}_i \to u^{(j)}_-$ in $C^2(\partial M)$ as $i \to \infty$. Then $u^{(j)}_- \in E_l(g_j(0)) = E_l(g) = E_k(g)$, $\norm{u^{(j)}_-}_{L^2(\partial M,g)} = 1$ and
\begin{align*}
Q_{(h_j, h_j\vert_{\partial M})}^{g,k}(u^{(j)}_-) = \lim_{i \to \infty} Q_{\omega_{ij}}^{g_j(t_i), l}(u^{(j)}_i) \geq 0.
\end{align*}
Again, taking a subsequence if necessary, $u^{(j)}_- \to u_-$ in $C^2(\partial M)$ as $j \to \infty$, and $u_- \in E_k(g)$, $\norm{u_-}_{L^2(\partial M,g)} = 1$. Moreover, 
\begin{align*}
Q_{(H, \eta)}^{g,k}(u_-) = \lim_{j \to \infty}Q_{(h_j, h_j\vert_{\partial M})}^{g,k}(u^{(j)}_-) \geq 0.
\end{align*}
Repeating the argument with $t \to 0^+$ produces $u_+ \in E_k(g)$, with $\norm{u_+}_{L^2(\partial M,g)} = 1$, such that $Q_{(H,\eta)}^{g,k}(u_+) \leq 0$. If $\dim E_k(g)\geq2$, the unit sphere of $E_k(g)$ is connected, and the conclusion follows from the intermediate value theorem. If $\dim E_k(g)=1$, then $u_+ = \pm u_-$ and, since $Q^{g,k}_{(H,\eta)}$ is quadratic, the two inequalities imply
\begin{align*}
Q^{g,k}_{(H,\eta)}(u_+)	= Q^{g,k}_{(H,\eta)}(u_-) =	0.
\end{align*}
\end{proof}

\subsection{Extremality under conformal metric deformations} \label{sec: Extremality under conformal metric deformations}

We next introduce the conformal counterparts of the Hilbert space $\mathcal H$, the quadratic forms $Q^{g,s}_{(H,\eta)}$, and the normalisation functional $N$. These objects allow us to formulate conformal extremality in parallel with the unrestricted setting. Although this section is not needed to prove the main results, we include it for completeness, since the corresponding arguments follow rather directly from the preceding analysis.

Given a Riemannian metric $g$ and a smooth function $\varphi$ on $M$, consider the Hilbert space
\begin{align*}
\mathcal{G}_\varphi = L^2(M, e^\varphi g) \times L^2(\partial M, e^\varphi g \vert_{\partial M})
\end{align*}
endowed with the product inner product. For a pair $(\psi, \xi) \in \mathcal{G}_\varphi$ and an integer $s \geq 1$, let $P_{(\psi, \xi)}^{e^\varphi g,s} : C^{\infty}(\partial M) \to \mathbb{R}$ be the quadratic form defined by
\begin{align*}
P_{(\psi, \xi)}^{e^\varphi g,s}(u) = \frac{(n-2)}{2} \int_M  \vert \nabla \hat u \vert^2 \psi \dV_{e^\varphi g} -  \frac{(n-1)}{2} \sigma_s(e^\varphi g) \int_{\partial M} u^2 \xi \dA_{e^\varphi g}.
\end{align*}
With this notation, the formula in Proposition \ref{derivative sigma conformal} can be written as
\begin{align*}
\dot \sigma_k(t_0) = P_{(\dot \varphi, \dot \varphi\vert_{\partial M})}^{e^{\varphi(t_0)}g, k}(u), 
\end{align*}
for any $u \in E_k(g(t_0))$ satisfying $\norm{u}_{L^2(\partial M,g(t_0))} = 1$.

We also denote by $N_g$ the restriction of the normalisation functional $N$ to the conformal class of $g$, that is:
\begin{align*}
N_g(\varphi) := N(e^\varphi g) = \Vol(M, e^\varphi g)^{\frac{2-n}{n}} \Vol(\partial M, e^\varphi g).
\end{align*}
The following lemma gives its first variation and identifies its gradient with respect to the Hilbert space $\mathcal{G}_\varphi$.

\begin{lem} \label{normalisation conformal}
Given $(-\varepsilon, \varepsilon) \ni t \mapsto \varphi(t)$ a smooth $1$-parameter family of smooth functions on $M$ with $\varphi(0) = \varphi$ and $\ddt \vert_{t = 0} \varphi(t) = \dot{\varphi}$, we have
\begin{align*}
\left. \ddt \right\vert_{t = 0} N_g(\varphi(t)) = \frac{N_g(\varphi)}{2} \left[ \frac{2 - n}{\Vol(M, e^\varphi g)} \int_M \dot{\varphi} \dV_{e^\varphi g} + \frac{n - 1}{\Vol(\partial M, e^\varphi g)} \int_{\partial M} \dot{\varphi} \dA_{e^\varphi g} \right].
\end{align*}
Alternatively, regarding $N_g$ as a function of the pair $(\varphi, \varphi\vert_{\partial M})$, we can write
\begin{align*}
\nabla N_g(\varphi, \varphi\vert_{\partial M}) = \frac{N_g(\varphi)}{2} \left( \frac{2 - n}{\Vol(M, e^\varphi g)}, \frac{n - 1}{\Vol(\partial M, e^\varphi g)} \right),
\end{align*}
so that
\begin{align*}
\left.\ddt\right|_{t=0}N_g(\varphi(t)) &= \left\langle\!\left\langle \nabla N_g(\varphi,\varphi|_{\partial M}), (\dot{\varphi},\dot{\varphi}|_{\partial M}) \right\rangle\!\right\rangle_{\mathcal{G_\varphi}}. 
 \end{align*} 
\end{lem}

Lemma \ref{normalisation conformal} describes the infinitesimal constraint imposed by fixing the value of $N_g$. Thus, the admissible directions at $\varphi$ are precisely the pairs $(\psi,\xi)\in\mathcal{G}_\varphi$ that are orthogonal to $\nabla N_g(\varphi,\varphi\vert_{\partial M})$. As remarked previously in the context of general metric variations, not every such pair is directly induced by a smooth conformal path, since a conformal variation produces the same function in the interior and on the boundary, namely a pair of the form $(f,f\vert_{\partial M})$, with $f \in C^{\infty}(M)$. The next lemma shows that this restriction is harmless at the $L^2$ level.

\begin{lem} 
The subspace
\begin{align*}
S = \{ (f, f\vert_{\partial M}) \in \mathcal{G}_\varphi : f \in C^{\infty}(M) \}
\end{align*}
is dense in $\mathcal{G}_\varphi$.
\end{lem}

We can now combine the conformal first-variation formula, the characterisation of the tangent space to the level set of $N_g$, and the preceding density lemma. Arguing as in the proof of Proposition \ref{isotropic}, we obtain the following conformal counterpart of that result.

\begin{prop} 
If a metric $g$ is $\overline{\sigma}_k$-conformally extremal, $k \geq 1$, then for any pair $(\psi, \xi) \in \mathcal{G}_0$ which satisfies
\begin{align*} 
\langle\!\langle \nabla N_g(0, 0\vert_{\partial M}), (\psi, \xi) \rangle\!\rangle_{\mathcal{G}_0} = 0
\end{align*}
there exists $u \in E_k(g)$ with $\norm{u}_{L^2(\partial M,g)} = 1$ such that $P_{(\psi, \xi)}^{g,k}(u) = 0$.
\end{prop}

\section{Characterisation of extremal metrics and proof of Theorem A} \label{sec: Characterisation of extremal metrics and proof of Theorem A}

In this section, we use the variational results of the preceding section to obtain a characterisation of extremal metrics in terms of Steklov eigenfunctions. We then establish the corresponding statement in the $G$-invariant setting and use it to prove Theorem A.

\subsection{Characterisation of extremal metrics} \label{sec: Characterisation of extremal metrics}

We begin with the unrestricted problem. The following result, adapted from the arguments of Fraser--Schoen and Karpukhin--M\'etras, shows that extremality forces the existence of a finite family of Steklov eigenfunctions satisfying pointwise identities in the interior and on the boundary.

\begin{thm}[\cite{fraser, karpukhin2022}, adapted] \label{thm:extremal steklov eigenfunctions}
Let $M$ be a compact connected smooth manifold with boundary $\partial M \neq \emptyset$ and $\dim M = n \geq 2$. If a Riemannian metric $g$ is $\overline{\sigma}_k$-extremal for some $k \geq 1$ and $n \geq 3$, then there exist Steklov eigenfunctions $u_1, \dots, u_m \in C^{\infty}(\partial M)$ associated with $\sigma_k(g)$ such that
\begin{itemize}
\item[(i)] $\displaystyle \sum_{i=1}^m u_i^2 = \frac{n}{\sigma_k(g) \Vol(\partial M, g)}$ on $\partial M$;
\item[(ii)] $\displaystyle \sum_{i=1}^m \dd \hat u_i \otimes \dd \hat u_i = \frac{g}{\Vol(M, g)}$ on $M$. 
\end{itemize} 
If the metric $g$ is $\overline{\sigma}_k$-extremal and $n = 2$, then there are eigenfunctions $u_1, \dots, u_m \in C^{\infty}(\partial M)$ associated with $\sigma_k(g)$ such that
\begin{itemize}
\item[(i)] $\displaystyle \sum_{i=1}^m u_i^2 = 1$ on $\partial M$;
\item[(ii)] $\displaystyle \sum_{i=1}^m \left( \dd \hat u_i \otimes \dd \hat u_i - \frac{1}{2} \vert \nabla \hat u_i \vert^2 g \right) = 0$ on $M$. 
\end{itemize} 

Conversely, if there is a collection of Steklov eigenfunctions of $(M^n,g)$, $n \geq 2$, associated with $\sigma_k(g)$, satisfying (i) and (ii) above and either $\sigma_{k-1}(g) < \sigma_k(g)$ or $\sigma_k(g) < \sigma_{k+1}(g)$, then $g$ is $\overline{\sigma}_k$-extremal.
\end{thm}

\begin{proof}[Proof of Theorem \ref{thm:extremal steklov eigenfunctions}]
Let $M^n$ be a compact smooth manifold with boundary and suppose $g$ is a smooth Riemannian metric which is extremal for $\overline{\sigma}_k$, $k \geq 1$. Let 
\begin{align*}
L(u) = \left( - \dd \hat u \otimes \dd \hat u + \frac{1}{2} \vert \nabla \hat u \vert^2 g,  -\frac{\sigma_k(g)}{2} u^2  g\vert_{\partial M} \right) \in \mathcal{H},
\end{align*}
so that\begin{align*}
Q_{(H, \eta)}^{g, k}(u) = \langle\!\langle L(u), (H, \eta) \rangle\!\rangle_{\mathcal{H}}.
\end{align*}

Define $K \subset \mathcal H$ to be the convex hull of
\begin{align*}
\{-L(u) : u \in E_k(g)\}.
\end{align*} 
We claim that $\nabla N := \nabla N(g, g\vert_{\partial M})$ belongs to $K$. If not, then since $K$ is a closed convex cone contained in a finite-dimensional subspace of $\mathcal{H}$, we can separate $\nabla N$ and $K$: there exists $\tilde{v} \in \mathcal{H}$ such that
\begin{align*}
\langle\!\langle \tilde{v}, \nabla N \rangle\!\rangle_{\mathcal{H}} < 0 \quad \text{ and } \quad \langle\!\langle \tilde{v}, - L(u) \rangle\!\rangle_{\mathcal{H}} \geq 0 \text{ for all } u \in E_k(g).
\end{align*}
Now set $G = (g, g\vert_{\partial M}) \in \mathcal{H}$ and define $v$ by
\begin{align*}
v = \tilde{v} - \frac{\langle\!\langle \tilde{v}, \nabla N \rangle\!\rangle_{\mathcal{H}}}{\langle\!\langle \nabla N, G \rangle\!\rangle_{\mathcal{H}}} \, G,
\end{align*}
so that $\langle\!\langle v, \nabla N \rangle\!\rangle_{\mathcal{H}} = 0$. Then Proposition \ref{isotropic} ensures the existence of $u \in E_k(g)$ with $\norm{u}_{L^2(\partial M, g)} = 1$ and $Q_{v}^{g, k}(u) = 0$. Note that $\langle\!\langle \nabla N, G \rangle\!\rangle_{\mathcal{H}} > 0$ and 
\begin{align*}
\langle\!\langle G, -L(u) \rangle\!\rangle_{\mathcal{H}} &= -\frac{(n-2)}{2} \int_M \vert \nabla \hat u \vert^2 \dV_{g} + \frac{(n-1)}{2} \sigma_k(g) \int_{\partial M} u^2 \dA_{g} \\
&= \frac{\sigma_k(g)}{2} \int_{\partial M} u^2 \dA_{g} >0.
\end{align*}
Thus,
\begin{align*}
0 = Q_{v}^{g, k}(u) = \langle\!\langle v, L(u) \rangle\!\rangle_{\mathcal{H}} = \underbrace{\langle\!\langle \tilde{v}, L(u) \rangle\!\rangle_{\mathcal{H}}}_{ \leq 0} + \underbrace{\frac{\langle\!\langle \tilde{v}, \nabla N \rangle\!\rangle_{\mathcal{H}}}{\langle\!\langle \nabla N, G \rangle\!\rangle_{\mathcal{H}}}}_{ < 0} \underbrace{\langle\!\langle G , -L(u) \rangle\!\rangle_{\mathcal{H}}}_{ > 0} < 0,
\end{align*}
which is a contradiction.

Since $K$ is a cone, $\frac{2}{N(g)} \nabla N \in K$. Thus, absorbing the coefficients into the eigenfunctions, there exist $u_1, \dots, u_m \in E_k(g)$ such that
\begin{align*}
\begin{cases}
\displaystyle \sum_{i=1}^m \left( \dd \hat u_i \otimes \dd \hat u_i - \frac{1}{2} \vert \nabla \hat u_i \vert^2 g \right) = \frac{(2-n)}{n} \frac{g}{\Vol(M,g)} \\
\displaystyle \sum_{i=1}^m u_i^2 = \frac{2}{\sigma_k(g) \Vol(\partial M, g)}.
\end{cases}
\end{align*}
Up to rescaling the $u_i$, this already implies the result for $n=2$. If $n \geq 3$, taking the trace of the first equation gives $\sum_{i = 1}^m \vert \nabla \hat u_i \vert^2 = \frac{2}{\Vol(M, g)}$. Inserting this back into the first equation yields $\sum_{i=1}^m \dd \hat u_i \otimes \dd \hat u_i = \frac{2g}{n \Vol(M, g)}$. The conclusion now follows by multiplying each $u_i$ by $\sqrt{n/2}$.

Conversely, suppose $u_1, \dots, u_m \in E_k(g)$ satisfy (i) and (ii). Let $(-\varepsilon, \varepsilon) \ni t \mapsto g(t)$ be a smooth $1$-parameter family of Riemannian metrics with $g(0) = g$, $\ddt g(t) \vert_{t = 0} = h$ and such that $N(g(t))$ is constant with respect to $t$ (we can always normalise $g(t)$ for this condition to be satisfied, as in the proof of Proposition \ref{isotropic}). Then, by Lemma \ref{normalisation},
\begin{align*}
	\frac{2 - n}{n \Vol(M, g)} \int_M \langle g, h \rangle \dV_g + \frac{1}{\Vol(\partial M, g)} \int_{\partial M} \langle g\vert_{\partial M}, h\vert_{\partial M} \rangle \dA_g = 0.
\end{align*}

Suppose first that $n \geq 3$. Using (i) and (ii), we obtain
\begin{align*}
	\sum_{i = 1}^m Q_{(h, h\vert_{\partial M})}^{g, k}(u_i) &= \int_M \left\langle \sum_{i = 1}^m \left( -\dd \hat u_i \otimes \dd \hat u_i + \frac{1}{2} \vert \nabla \hat u_i \vert^2 g \right), h \right\rangle \dV_g \\
	&\quad - \frac{\sigma_k(g)}{2} \int_{\partial M} \sum_{i = 1}^m u_i^2 \langle g\vert_{\partial M}, h\vert_{\partial M} \rangle \dA_g \\
	&= \frac{n - 2}{2 \Vol(M, g)} \int_M \langle g, h \rangle \dV_g - \frac{n}{2 \Vol(\partial M, g)} \int_{\partial M} \langle g\vert_{\partial M}, h\vert_{\partial M} \rangle \dA_g \\
	&= 0.
\end{align*}
For $n = 2$, the interior term vanishes by (ii), while (i) gives $\sum_{i = 1}^m u_i^2 = 1$ on $\partial M$. Since $N(g(t)) = \Vol(\partial M, g(t))$ is constant, the boundary term also vanishes. Hence, in either dimension,
\begin{align*}
\sum_{i = 1}^m Q_{(h, h\vert_{\partial M})}^{g, k}(u_i) = 0.
\end{align*}

Letting $E = \operatorname{span} \{u_1, \dots, u_m\} \subset E_k(g)$, the preceding identity implies that there exist $u_{\pm} \in E$, with $\norm{u_{\pm}}_{L^2(\partial M, g)} = 1$, such that $\pm Q_{(h, h\vert_{\partial M})}^{g, k}(u_{\pm}) \leq 0$.

Suppose now that $\sigma_{k - 1}(g) < \sigma_k(g)$. Let $P_t$ denote the spectral projection onto the first $k$ Steklov eigenfunctions of $g(t)$, and set $u_t = u_+ - P_t u_+$. As in the proof of Proposition \ref{derivative sigma}, the spectral gap ensures that $P_t$ varies smoothly, while the min--max characterisation gives
\begin{align*}
\sigma_k(g(t)) \leq \frac{\int_M \vert \nabla \hat u_t \vert_{g(t)}^2 \, \dV_{g(t)}}{\int_{\partial M} u_t^2 \, \dA_{g(t)}},
\end{align*}
with equality at $t = 0$. Differentiating the quotient at $t = 0$ as in Proposition \ref{derivative sigma}, we obtain
\begin{align*}
\limsup_{t \to 0^+} \frac{\sigma_k(g(t)) - \sigma_k(g)}{t} \leq Q_{(h, h\vert_{\partial M})}^{g,k}(u_+) \leq 0.
\end{align*}
Similarly, using $u_-$ instead of $u_+$ and letting $t \to 0^-$, we obtain
\begin{align*}
\liminf_{t \to 0^-} \frac{\sigma_k(g(t)) - \sigma_k(g)}{t} \geq Q_{(h, h\vert_{\partial M})}^{g,k}(u_-) \geq 0.
\end{align*}
Consequently, $\overline{\sigma}_k(g(t)) \leq \overline{\sigma}_k(g) + o(t)$ as $t \to 0$. 

The case $\sigma_k(g) < \sigma_{k + 1}(g)$ is treated similarly, using the spectral projection onto the first $k + 1$ eigenspaces and obtaining the reverse extremality inequality. Thus, in either case, $g$ is $\overline{\sigma}_k$-extremal.
\end{proof}

We now apply the same argument to $G$-invariant metric deformations. The variational framework developed in the previous section yields the following $G$-invariant counterpart.

\begin{thm} \label{thm:G-extremal steklov eigenfunctions}
Let $M^n$ be a compact connected smooth manifold with boundary $\partial M \neq \emptyset$ and $\dim M = n \geq 2$. Suppose $G$ is a compact Lie group acting smoothly and effectively on $M$, and assume the action is boundary-nontransitive. If a $G$-invariant Riemannian metric $g$ is $\overline{\sigma}_k^G$-extremal for some $k \geq 1$ and $n \geq 3$, then there exist $G$-invariant Steklov eigenfunctions $u_1, \dots, u_m \in C^{\infty}(\partial M)^G$ associated with $\sigma_k^G(g)$ such that
\begin{itemize}
\item[(i)] $\displaystyle \sum_{i=1}^m u_i^2 = \frac{n}{\sigma_k^G(g) \Vol(\partial M, g)}$ on $\partial M$;
\item[(ii)] $\displaystyle \sum_{i=1}^m \dd \hat u_i \otimes \dd \hat u_i = \frac{g}{\Vol(M, g)}$ on $M$. 
\end{itemize} 
If the metric $g$ is $\overline{\sigma}_k^G$-extremal and $n = 2$, then there are $G$-invariant eigenfunctions $u_1, \dots, u_m \in C^{\infty}(\partial M)^G$ associated with $\sigma_k^G(g)$ such that
\begin{itemize}
\item[(i)] $\displaystyle \sum_{i=1}^m u_i^2 = 1$ on $\partial M$;
\item[(ii)] $\displaystyle \sum_{i=1}^m \left( \dd \hat u_i \otimes \dd \hat u_i - \frac{1}{2} \vert \nabla \hat u_i \vert^2 g \right) = 0$ on $M$.
\end{itemize} 
\end{thm}

\begin{proof}
The proof follows the same argument as that of Theorem \ref{thm:extremal steklov eigenfunctions}, with a few standard observations concerning $G$-invariance.
	
Let $\mathcal H^G$ denote the closed subspace of $\mathcal H$ consisting of $G$-invariant tensor pairs, and let $E_k^G(g)$ be the eigenspace associated with $\sigma_k^G(g)$. Since harmonic extension commutes with the action of $G$, every $u \in E_k^G(g)$ has a $G$-invariant harmonic extension. Consequently, the tensor pair
\begin{align*}
L_G(u) = \left(-\dd \hat u \otimes \dd \hat u + \frac{1}{2} \vert \nabla \hat u \vert^2 g, -\frac{\sigma_k^G(g)}{2}u^2 g\vert_{\partial M}\right)
\end{align*}
belongs to $\mathcal H^G$. The gradient $\nabla N(g, g\vert_{\partial M})$ is also $G$-invariant.

We first observe that Proposition \ref{isotropic} remains valid when $\mathcal H$ and $E_k(g)$ are replaced by $\mathcal H^G$ and $E_k^G(g)$. Indeed, the density argument of Lemma \ref{denseness} carries over by averaging smooth approximations over $G$, and the scalar correction used in the projection preserves $G$-invariance. Moreover, the first-variation argument of Proposition \ref{derivative sigma} applies to smooth families of $G$-invariant metrics almost verbatim. Thus, for every $V \in \mathcal H^G$ orthogonal to $\nabla N(g, g\vert_{\partial M})$, there exists $u \in E_k^G(g)$ with $\norm{u}_{L^2(\partial M, g)} = 1$ such that $\langle\!\langle L_G(u), V \rangle\!\rangle_{\mathcal H} = 0$.
	
Now let $K^G$ be the closed convex cone generated by $\{-L_G(u) : u \in E_k^G(g)\}$ in $\mathcal H^G$. The separation argument in the proof of Theorem \ref{thm:extremal steklov eigenfunctions} applies without modification and shows that $\nabla N(g, g\vert_{\partial M}) \in K^G$. Consequently, there exist $G$-invariant eigenfunctions $u_1, \dots, u_m \in E_k^G(g)$ such that
\begin{align*}
\sum_{i = 1}^m \left(\dd \hat u_i \otimes \dd \hat u_i - \frac{1}{2} \vert \nabla \hat u_i \vert^2 g\right) = \frac{(2 - n)}{n}\frac{g}{\Vol(M, g)}, \qquad \sum_{i = 1}^m u_i^2 = \frac{2}{\sigma_k^G(g)\Vol(\partial M, g)}.
\end{align*}
Taking the trace of the first identity when $n \geq 3$, and rescaling the eigenfunctions as in the proof of Theorem \ref{thm:extremal steklov eigenfunctions}, gives the desired identities. The case $n = 2$ follows directly by rescaling the boundary identity.
\end{proof}

\subsection{The proof of Theorem \ref{thm: no extremal metrics}} \label{sec: The proof of Theorem no extremal metrics}

We now prove Theorem \ref{thm: no extremal metrics} as a direct consequence of the preceding characterisation of $G$-invariant extremal metrics.

Under the hypotheses of the theorem, suppose that $g$ is a smooth $G$-invariant $\overline{\sigma}_k^G$-extremal metric on $M$, for some $k\geq1$. Since $G$ has positive dimension and acts effectively on $M$, there exists a point $p\in \interior M$ whose orbit has positive dimension. We may therefore choose a nonzero vector $X_p\in T_p(G\cdot p)$ tangent to the $G$-orbit at $p$.

By Theorem \ref{thm:G-extremal steklov eigenfunctions}, there exist $G$-invariant Steklov eigenfunctions $u_1,\ldots,u_m$ associated with $\sigma_k^G(g)$ such that 
\begin{align} \label{eq: G eigenfunction equation n=3}
\sum_{i=1}^m \dd\hat u_i \otimes \dd\hat u_i =\frac{g}{\Vol(M,g)}.
\end{align}
Since $g$ is $G$-invariant, the harmonic extension $\hat u_i$ of each $u_i$ is $G$-invariant. Hence $\dd(\hat u_i)_p X_p = 0$ for every $i=1,\ldots,m$. Evaluating \eqref{eq: G eigenfunction equation n=3} on the pair $(X_p,X_p)$ gives
\begin{align*}
0 = \frac{g(X_p,X_p)}{\Vol(M,g)},
\end{align*}
which is impossible since $X_p \neq 0$.

\section{The proof of Theorem \ref{thm: infinite sup}} \label{sec: The proof of Theorem infinite sup}

The proof proceeds by constructing a local model near a principal boundary orbit of positive dimension. On shrinking $G$-invariant neighborhoods of this orbit, the normalised $G$-invariant Steklov--Neumann eigenvalues diverge. After comparing the local geometry with the model geometry, we use a conformal collapse of the complement to transfer this divergence to the Steklov spectrum of $M$. For clarity, we divide the proof into five steps.

\bigskip

\textbf{Step 1: The local model}

\bigskip

In $\mathbb R^m$, $m \geq 2$, consider the closed half-ball
\begin{align*}
	\mathbb{B}^+ = \{x \in \mathbb R^{m - 1} \times [0, +\infty) : \vert x \vert \leq 1\},
\end{align*}
and, for $\rho > 0$, set
\begin{align*}
	\mathbb{B}_\rho^+ = \rho \mathbb{B}^+ = \{x \in \mathbb R^{m - 1} \times [0, +\infty) : \vert x \vert \leq \rho\}.
\end{align*}
We distinguish the flat and spherical portions of the boundary of $\mathbb{B}^+$, defined by
\begin{align*}
	\partial_V \mathbb{B}^+ &= \{(x', 0) \in \mathbb R^{m - 1} \times \{0\} : \vert x' \vert < 1\}, \\
	\partial_L \mathbb{B}^+ &= \{x \in \mathbb R^{m - 1} \times (0, +\infty) : \vert x \vert = 1\},
\end{align*}
and analogously for $\mathbb{B}_\rho^+$. Their closures cover $\partial \mathbb{B}^+$ and intersect along the equatorial $(m - 2)$-sphere.

As in Section \ref{sec: Steklov-Neumann problem}, let $\mu_k(\mathbb{B}_\rho^+)$ denote the $k$-th eigenvalue of the Steklov--Neumann problem
\begin{align*}
\begin{cases}
\Delta u =0 & \text{in } \mathbb{B}_\rho^+, \\
\partial_\nu u = \mu u & \text{on }\partial_V \mathbb{B}_\rho^+, \\
\partial_\nu u = 0 & \text{on } \partial_L \mathbb{B}_\rho^+.
\end{cases}
\end{align*}

\begin{lem}
For any $\rho > 0$ and $k \geq 0$,
\begin{align*}
\mu_k(\mathbb{B}_\rho^+) = \frac{1}{\rho} \mu_k(\mathbb{B}^+).
\end{align*}
In particular, for any $k \geq 1$
\begin{align*}
\lim_{\rho \to 0} \mu_k(\mathbb{B}_\rho^+) = + \infty.
\end{align*}
\end{lem}

\begin{proof}
The Steklov--Neumann eigenvalues $\mu_k(\mathbb{B}^+)$ admit the variational characterisation
\begin{align} \label{steklov-neumann characterisation}
\mu_k(\mathbb{B}^+) = \min_{V \in \mathcal{V}_{k+1}(\mathbb{B}^+)} \max_{0 \neq u \in V} \frac{\displaystyle \int_{\mathbb{B}^+} \vert \nabla u \vert^2 \, \dd x}{\displaystyle \int_{\partial_V \mathbb{B}^+} u^2 \, \dd x'}, 
\end{align}
where $\mathcal{V}_k(\mathbb{B}^+)$ denotes the set of all $k$-planes in $H^1(\mathbb{B}^+)$ and $\dd x'$ is the Lebesgue measure in $\mathbb{R}^{m-1} \times \{0\} \subset \mathbb{R}^m$.

Consider the homothety $h_\rho: \mathbb{B}^+ \to \mathbb{B}_\rho^+$ given by $h_\rho(y)=\rho y$. This map induces an isomorphism $h_\rho^\ast : H^1(\mathbb{B}_\rho^+) \to H^1(\mathbb{B}^+)$ via $h_\rho^\ast(u) = u \circ h_\rho$. Set $v = u \circ h_\rho$.  Then
\begin{align*}
\nabla v(y) = \rho \nabla u(\rho y),
\end{align*}
and therefore, the change of variable formula gives:
\begin{align*}
\int_{\mathbb{B}_\rho^+} \vert \nabla u \vert^2 \, \dd x = \rho^{m-2} \int_{\mathbb{B}^+} \vert \nabla v \vert^2 \, \dd y.
\end{align*}
Analogously,
\begin{align*}
\int_{\partial_V \mathbb{B}_\rho^+} u^2 \, \dd x' = \rho^{m-1}\int_{\partial_V \mathbb{B}^+} v^2 \, \dd y'.
\end{align*}
Hence
\begin{align*}
\frac{\displaystyle\int_{\mathbb{B}_\rho^+} \vert \nabla u \vert^2 \, \dd x}{\displaystyle\int_{\partial_V \mathbb{B}_\rho^+} u^2 \, \dd x'} = \frac{1}{\rho} \frac{\displaystyle\int_{\mathbb{B}^+} \vert \nabla v \vert^2 \, \dd y}{\displaystyle\int_{\partial_V \mathbb{B}^+} v^2 \, \dd y'}.
\end{align*}
Moreover, $h_\rho^\ast$ induces a bijection between $\mathcal{V}_{k+1}(\mathbb{B}_\rho^+)$ and $\mathcal{V}_{k+1}(\mathbb{B}^+)$. Therefore, applying \eqref{steklov-neumann characterisation} on both domains yields
\begin{align*}
\mu_k(\mathbb{B}_\rho^+) = \frac{1}{\rho} \mu_k(\mathbb{B}^+).
\end{align*}
Since $\mathbb{B}^+$ is connected, $\mu_0(\mathbb{B}^+) = 0$ is simple, and hence $\mu_k(\mathbb{B}^+) > 0$ for every $k \geq 1$. The last assertion follows.
\end{proof}

\bigskip

\textbf{Step 2: The product model}

\bigskip

Let $H$ be a closed subgroup of $G$, set $q = \dim(G/H)$, and let $m \geq 2$. We put $n = q + m$. Endow $G/H$ with a smooth $G$-invariant Riemannian metric $h$ and, for $\rho > 0$, consider
\begin{align*}
\Omega_\rho = (G/H) \times \mathbb{B}_\rho^+
\end{align*}
endowed with the product metric $g_\rho = h \oplus g_{\mathrm{euc}}$. We let $G$ act on $\Omega_\rho$ by left multiplication on the first factor and trivially on the second, that is,
\begin{align*}
\gamma \cdot (y, x) = (\gamma y, x), \qquad (y, x) \in (G/H) \times \mathbb{B}_\rho^+
\end{align*}
The boundary of $\Omega_\rho$ is decomposed as $\partial \Omega_\rho = \overline{\partial_V \Omega_\rho} \cup \overline{\partial_L \Omega_\rho}$, where
\begin{align*}
\partial_V \Omega_\rho &= (G/H) \times \partial_V \mathbb{B}_\rho^+,\\
\partial_L \Omega_\rho &= (G/H) \times \partial_L \mathbb{B}_\rho^+.
\end{align*}

As in Section \ref{sec: Steklov-Neumann problem}, let $\mu_k^G(\Omega_\rho)$ denote the $k$-th $G$-invariant eigenvalue of the Steklov--Neumann problem
\begin{align*}
\begin{cases}
\Delta_{g_\rho} u = 0 & \text{in } \Omega_\rho,\\
\partial_\nu u = \mu u & \text{on } \partial_V \Omega_\rho,\\
\partial_\nu u = 0 & \text{on } \partial_L \Omega_\rho.
\end{cases}
\end{align*}

The following lemma shows that these eigenvalues are easily computable, since the action of $G$ on $G/H$ is transitive.

\begin{lem}
For every $\rho > 0$ and $k \geq 0$
\begin{align*}
\mu_k^G(\Omega_\rho) = \mu_k(\mathbb{B}_\rho^+) = \frac{1}{\rho}\mu_k(\mathbb{B}^+).
\end{align*}
Moreover,
\begin{align*}
\mu_k^G(\Omega_\rho)\Vol(\Omega_\rho,g_\rho)^{\frac{2-n}{n}}\Vol(\partial_V\Omega_\rho,g_\rho) = C \rho^{-\frac{2q}{n}},
\end{align*}
where $C = C(m, q, k, h) \geq 0$ is independent of $\rho$ and $C > 0$ for $k \geq 1$. In particular, if $q > 0$ and $k \geq 1$, then
\begin{align*}
\lim_{\rho \to 0} \mu_k^G(\Omega_\rho) \Vol(\Omega_\rho, g_\rho)^{\frac{2 - n}{n}} \Vol(\partial_V \Omega_\rho, g_\rho) = +\infty.
\end{align*}
\end{lem}

\begin{proof}
Since $G$ acts transitively on $G/H$, every $G$-invariant function $u\in H^1(\Omega_\rho)^G$ is of the form $u(y,x) = v(x)$ for some unique $v \in H^1(\mathbb{B}_\rho^+)$. Conversely, every function of this form is $G$-invariant. Hence the corresponding map
\begin{align*}
H^1(\mathbb{B}_\rho^+) \ni v \mapsto u \in H^1(\Omega_\rho)^G
\end{align*}
is an isomorphism.
	
Let $\alpha_h = \Vol(G/H,h)$. For $u(y,x)=v(x)$, the product structure gives
\begin{align*}
\int_{\Omega_\rho} \vert \nabla u \vert^2 \, \dV_{g_\rho} = \alpha_h \int_{\mathbb{B}_\rho^+} \vert \nabla v \vert^2 \, \dd x
\end{align*}
and
\begin{align*}
\int_{\partial_V \Omega_\rho} u^2 \, \dA_{g_\rho} = \alpha_h \int_{\partial_V \mathbb{B}_\rho^+} v^2 \, \dd x'.
\end{align*}
Therefore,
\begin{align*}
\frac{\displaystyle \int_{\Omega_\rho} \vert \nabla u \vert^2 \, \dV_{g_\rho}}{\displaystyle \int_{\partial_V \Omega_\rho} u^2 \, \dA_{g_\rho}} = \frac{\displaystyle \int_{\mathbb{B}_\rho^+} \vert \nabla v \vert^2 \, \dd x}{\displaystyle \int_{\partial_V \mathbb{B}_\rho^+} v^2 \, \dd x'}.
\end{align*}
The preceding isomorphism induces a bijection between $(k+1)$-planes in $H^1(\mathbb{B}_\rho^+)$ and $(k+1)$-planes in $H^1(\Omega_\rho)^G$. The variational characterisations of the two Steklov--Neumann spectra and Step 1 therefore yield
\begin{align*}
\mu_k^G(\Omega_\rho) = \mu_k(\mathbb{B}_\rho^+) = \frac{1}{\rho} \mu_k(\mathbb{B}^+).
\end{align*}
	
We now compute the normalised $G$-invariant eigenvalue $\mu_k^G(\Omega_\rho)$ for $k \geq 1$. Since
\begin{align*}
\Vol(\mathbb{B}_\rho^+) = \rho^m \Vol(\mathbb{B}^+)
\end{align*}
and
\begin{align*}
\Vol(\partial_V \mathbb{B}_\rho^+) = \rho^{m - 1} \Vol(\partial_V \mathbb{B}^+),
\end{align*}
we have
\begin{align*}
\Vol(\Omega_\rho, g_\rho) = \alpha_h \rho^m \Vol(\mathbb{B}^+)
\end{align*}
and
\begin{align*}
\Vol(\partial_V \Omega_\rho, g_\rho) = \alpha_h \rho^{m - 1} \Vol(\partial_V \mathbb{B}^+).
\end{align*}
Consequently,
\begin{align*}
\mu_k^G(\Omega_\rho) \Vol(\Omega_\rho, g_\rho)^{\frac{2 - n}{n}} \Vol(\partial_V \Omega_\rho,g_\rho) = C \rho^{-\frac{2q}{n}},
\end{align*}
where
\begin{align*}
C = \mu_k(\mathbb{B}^+) \left( \alpha_h \Vol(\mathbb{B}^+) \right)^{\frac{2 - n}{n}} \alpha_h \Vol(\partial_V \mathbb{B}^+) > 0.
\end{align*}
The conclusion follows whenever $q > 0$ and $k \geq 1$.
\end{proof}

\textbf{Step 3: Quasi-isometric metrics}

\bigskip

We next show that the preceding estimates are stable under quasi-isometric perturbations of the metric. Let $\Omega$ be a compact Lipschitz manifold of dimension $n \geq 3$, whose boundary is decomposed as
\begin{align*}
	\partial \Omega = \overline{\partial_V \Omega} \cup \overline{\partial_L \Omega},
\end{align*}
where $\partial_V \Omega$ and $\partial_L \Omega$ are disjoint relatively open smooth subsets of $\partial\Omega$ whose closures meet along a smooth codimension-two submanifold. We assume that every connected component of $\Omega$ meets $\partial_V \Omega$, as in Remark \ref{rem: disconnected problem}, and that a compact Lie group $G$ acts smoothly on $\Omega$, preserving both boundary portions.

As in Section \ref{sec: Steklov-Neumann problem}, given a $G$-invariant Riemannian metric $g$ on $\Omega$, let $\mu_k^G(\Omega, g)$ denote the $k$-th $G$-invariant eigenvalue of the Steklov--Neumann problem
\begin{align*}
\begin{cases}
\Delta_g u = 0 & \text{in } \Omega,\\
\partial_\nu u = \mu u & \text{on } \partial_V \Omega,\\
\partial_\nu u = 0 & \text{on } \partial_L \Omega.
\end{cases}
\end{align*}
As remarked in that section, in general there may be only finitely many such eigenvalues if no further assumptions on the $G$-action are imposed. For the purely variational comparison below, however, no such additional assumptions are needed.

\begin{lem} \label{quasi-isometry steklov-neumann}
Let $g_1$ and $g_2$ be two $G$-invariant Riemannian metrics on $\Omega$	such that
\begin{align*}
C^{-1} g_1 \leq g_2 \leq C g_1
\end{align*}
for some $C \geq 1$. Then, for every $k \geq 0$ for which the corresponding $G$-invariant Steklov--Neumann eigenvalue is defined and finite, we have
\begin{align*}
C^{-\frac{2n+1}{2}} \mu_k^G(\Omega,g_1)	\leq \mu_k^G(\Omega,g_2) \leq C^{\frac{2n+1}{2}} \mu_k^G(\Omega,g_1).
\end{align*}
Moreover, letting
\begin{align*}
\overline{\mu}_k^G(\Omega,g) = \mu_k^G(\Omega,g) \Vol(\Omega,g)^{\frac{2-n}{n}} 		\Vol(\partial_V\Omega,g),
\end{align*}
then
\begin{align*}
C^{-(2n-1)}\overline{\mu}_k^G(\Omega,g_1) \leq \overline{\mu}_k^G(\Omega,g_2) \leq 		C^{2n-1}\overline{\mu}_k^G(\Omega,g_1).
\end{align*}
\end{lem}

\begin{proof}
The assumption implies
\begin{align*}
C^{-1}\vert \nabla u \vert_{g_1}^2 \leq \vert \nabla u \vert_{g_2}^2 \leq C\vert \nabla u \vert_{g_1}^2
\end{align*}
for every $u\in H^1(\Omega)$. The corresponding volume elements satisfy
\begin{align*}
C^{-\frac{n}{2}} \dV_{g_1} \leq \dV_{g_2} \leq C^{\frac{n}{2}} \dV_{g_1},
\end{align*}
while, on $\partial_V\Omega$,
\begin{align*}
C^{-\frac{n-1}{2}} \dA_{g_1} \leq \dA_{g_2} \leq C^{\frac{n-1}{2}} \dA_{g_1}.
\end{align*}
It follows that, for every nonzero $u\in H^1(\Omega)^G$,
\begin{align*}
C^{-\frac{2n+1}{2}}	\frac{\displaystyle\int_\Omega \vert \nabla u\vert_{g_1}^2 \, \dV_{g_1}}		{\displaystyle\int_{\partial_V\Omega}u^2 \, \dA_{g_1}}
\leq
\frac{\displaystyle\int_\Omega \vert \nabla u\vert_{g_2}^2 \, \dV_{g_2}}		{\displaystyle\int_{\partial_V\Omega}u^2 \, \dA_{g_2}}
\leq
C^{\frac{2n+1}{2}} \frac{\displaystyle\int_\Omega \vert \nabla u\vert_{g_1}^2 \, \dV_{g_1}}		{\displaystyle\int_{\partial_V\Omega}u^2 \, \dA_{g_1}}.
\end{align*}
The quasi-isometry assumption implies that the Sobolev spaces associated with $g_1$ and $g_2$ coincide as vector spaces and have equivalent norms. Moreover, $G$-invariance is independent of the metric. Hence the same $(k+1)$-dimensional $G$-invariant subspaces of $H^1(\Omega)^G$ may be used in the variational characterisations for both metrics. We therefore obtain
\begin{align*}
C^{-\frac{(2n+1)}{2}} \mu_k^G(\Omega,g_1)	\leq \mu_k^G(\Omega,g_2) \leq C^{\frac{2n+1}{2}} \mu_k^G(\Omega,g_1).
\end{align*}
The estimates for the interior and boundary volumes follow directly	from the corresponding estimates for the volume elements. Indeed,
\begin{align*}
C^{-\frac{n}{2}}\Vol(\Omega,g_1) \leq \Vol(\Omega,g_2) \leq C^{\frac{n}{2}}\Vol(\Omega,g_1),
\end{align*}
and 
\begin{align*}
C^{-\frac{(n-1)}{2}}\Vol(\partial_V\Omega,g_1) \leq \Vol(\partial_V\Omega,g_2) \leq C^{\frac{n-1}{2}}\Vol(\partial_V\Omega,g_1).
\end{align*}
Since the exponent $(2-n)/n$ is nonpositive, combining the preceding three estimates yields the stated result.
\end{proof}

\textbf{Step 4: A neighbourhood of a boundary orbit}

\bigskip

We now show that the product model of Step 2 describes, up to a quasi-isometric error tending to zero, a sufficiently small neighbourhood of a suitable boundary orbit in $M$.

Let $N$ be a connected component of $\partial M$ on which the action of $G_N$ is nontransitive. Since $G$ permutes the finitely many connected components of $\partial M$, the subgroup $G_N$ has finite index in $G$, and hence $\dim G_N = \dim G = d \geq 1$.

Choose $p \in N$ in the principal stratum of the $G_N$-action on $N$, and set
\begin{align*}
H = G_p, \qquad q = \dim(G \cdot p), \qquad m = n - q.
\end{align*}
Every element of $G$ fixing $p$ preserves $N$, so $H \subseteq G_N$. Moreover, since $G_N$ has finite index in $G$, we have $T_p (G \cdot p) = T_p (G_N \cdot p)$. Hence $q = \dim(G_N \cdot p)$, and the nontransitivity of the $G_N$-action on the $(n-1)$-dimensional manifold $N$ implies that $q \leq n - 2$.

We also have $q > 0$. Indeed, otherwise the principal $G_N$-orbits would be zero-dimensional. Since $G^0 \subseteq G_N$, the identity component $G^0$ would fix the principal stratum of $N$ pointwise, and hence all of $N$. Thus, for every $\gamma \in G^0$ and $x \in N$, the differential $\dd \gamma_x$ is the identity on $T_x N$. Since $\gamma$ also preserves the inward unit normal to $\partial M$, it follows that $\dd \gamma_x = \id_{T_x M}$. An isometry of the connected manifold $(M, g_0)$ is determined by its value and differential at one point, so $\gamma = \id_M$, contradicting effectiveness and $\dim G > 0$. Consequently, $1 \leq q \leq n - 2$, and therefore $m \geq 2$.

Let
\begin{align*}
V_p = T_p (G \cdot p)^\perp \cap T_p N,
\end{align*}
where orthogonality is taken with respect to $g_0$, and let $\eta_p$ be the inward unit normal to $\partial M$ at $p$. Then
\begin{align*}
T_p M = T_p (G \cdot p) \oplus V_p \oplus \mathbb{R}\eta_p
\end{align*}
is an orthogonal decomposition, with $\dim V_p = m - 1 \geq 1$.

Since $p$ belongs to the principal stratum of the $G_N$-action on $N$, the slice representation of $H = G_p$ on $V_p$ is trivial. Moreover, the action of $H$ on $\mathbb{R}\eta_p$ is trivial, since every element of $H$ preserves the inward normal half-line. Hence $H$ acts trivially on $V_p \oplus \mathbb{R} \eta_p$.

We endow $V_p \oplus \mathbb R\eta_p$ with the Euclidean metric $g_{\mathrm{euc}}$ induced by $g_0(p)$. Choosing an isometry from this space onto $\mathbb R^m$ that sends $\eta_p$ to the last coordinate vector, we identify
\begin{align*}
\mathbb{B}_\delta^+ = \left\{v + t\eta_p \in V_p \oplus \mathbb R\eta_p : \vert v \vert^2 + t^2 \leq \delta^2,\ t \geq 0\right\}
\end{align*}
with the standard Euclidean half-ball of radius $\delta$.

\begin{figure}[h]
\centering
\includegraphics[width = 8cm]{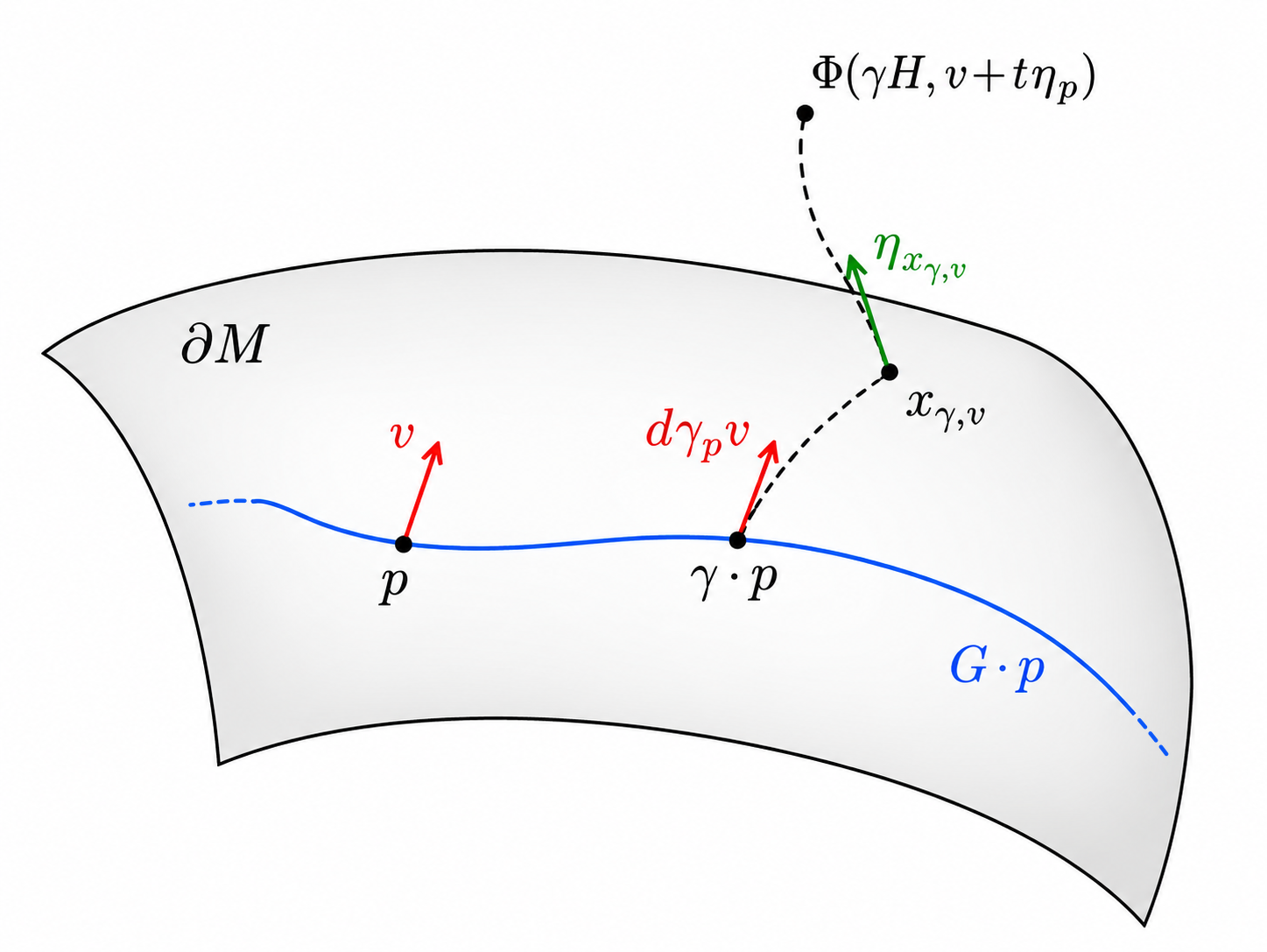}
\caption{Construction of the $G$-equivariant map $\Phi$ near the boundary orbit $G \cdot p$. The black dashed curves represent the geodesic in $\partial M$ from $\gamma\cdot p$ to $x_{\gamma, v} = \exp^{\partial M}_{\gamma \cdot p}(d \gamma_p v)$, followed by the inward normal geodesic in $M$ joining $x_{\gamma, v}$ to $\Phi(\gamma H, v + t \eta_p)$.}
\label{construction}
\end{figure}

Let $\eta$ denote the inward unit normal field along $\partial M$. For sufficiently small $\delta > 0$, define
\begin{align*}
\Phi : (G/H) \times \mathbb B_\delta^+ &\longrightarrow M,\\
(\gamma H, v + t\eta_p) &\longmapsto \exp_{x_{\gamma,v}}^M(t\eta_{x_{\gamma,v}}),
\end{align*}
where
\begin{align*}
x_{\gamma, v} = \exp_{\gamma \cdot p}^{\partial M}(\dd \gamma_p v).
\end{align*}
Here $\exp^{\partial M}$ is the exponential map of the metric induced by $g_0$ on $\partial M$ and $\exp^M$ is the exponential map of $g_0$ on $M$. This construction is illustrated in Figure \ref{construction}. 

Since $H$ acts trivially on $V_p$ and fixes $\eta_p$, the definition is independent of the representative $\gamma$ of $\gamma H$. Moreover, $\Phi$ is $G$-equivariant for the action
\begin{align*}
\gamma \cdot(y, x) = (\gamma y, x), \qquad (y, x) \in(G/H) \times \mathbb B_\delta^+.
\end{align*}

By construction, $\Phi(\gamma H,0) = \gamma\cdot p$, and the differential of $\Phi$ at $(eH,0)$ identifies the decomposition
\begin{align*}
T_{eH}(G/H)\oplus V_p \oplus \mathbb{R} \eta_p
\end{align*}
with the orthogonal decomposition
\begin{align*}
T_pM = T_p(G \cdot p) \oplus V_p \oplus \mathbb R\eta_p.
\end{align*}
In particular, this differential is invertible. By equivariance, the same holds at every point of $(G/H)\times\{0\}$. 
Since $\Phi$ is injective on the zero section and $G/H$ is compact, the inverse function theorem implies that, after decreasing $\delta$ if necessary, $\Phi$ is a smooth $G$-equivariant embedding whose image contains a neighbourhood of $G \cdot p$ in $M$. Moreover,
\begin{align*}
\Phi^{-1}(\partial M) = (G/H)\times\left(\mathbb B_\delta^+\cap(V_p\oplus\{0\})\right).
\end{align*}

We next compare $\Phi^\ast g_0$ with the product metric arising from the geometry of the orbit at $p$. Let $h$ be the $G$-invariant metric on $G/H$ for which the orbit map $\gamma H\mapsto\gamma\cdot p$ is an isometry onto $G\cdot p$, endowed with the metric induced by $g_0$. We set
\begin{align*}
g_{\mathrm{mod}} = h \oplus g_{\mathrm{euc}}.
\end{align*}
By the construction of $\Phi$, its differential at $(eH, 0)$ is an isometry with respect to $g_{\mathrm{mod}}$ and $g_0$. Thus $\Phi^\ast g_0 = g_{\mathrm{mod}}$ at $(eH, 0)$. Since $\Phi$ is $G$-equivariant and both metrics are $G$-invariant, it follows that $\Phi^\ast g_0 = g_{\mathrm{mod}}$ on $(G/H) \times \{0\}$.

For each $0 < \rho < \delta$, define
\begin{align*}
\varepsilon(\rho) = \sup \left\{ \left\vert (\Phi^\ast g_0 - g_{\mathrm{mod}})_{(y, x)}(X, X) \right\vert :
\begin{array}{l}
y \in G/H,\  x \in \mathbb{B}_\rho^+, \\
X \in T_y(G/H) \oplus \mathbb{R}^m,\ 
(g_{\mathrm{mod}})_{(y, x)}(X, X) = 1
\end{array}
\right\}.
\end{align*}
Since the two metrics agree on $(G/H) \times \{0\}$, compactness of $G/H$ and continuity of the metrics imply that $\varepsilon(\rho) \to 0$ as $\rho \to 0$. Moreover, the definition of $\varepsilon$ gives
\begin{align*}
(1 - \varepsilon(\rho))g_{\mathrm{mod}} \leq \Phi^\ast g_0 \leq (1 + \varepsilon(\rho))g_{\mathrm{mod}}
\end{align*}
on $(G/H) \times \mathbb{B}_\rho^+$. For $\rho$ sufficiently small, $\varepsilon(\rho) < 1$, so letting
\begin{align*}
C_\rho = \frac{1}{1 - \varepsilon(\rho)},
\end{align*}
we have $C_\rho \to 1$ as $\rho \to 0$ and
\begin{align*}
C_\rho^{-1}g_{\mathrm{mod}} \leq \Phi^\ast g_0 \leq C_\rho g_{\mathrm{mod}}
\end{align*}
on $(G/H) \times \mathbb{B}_\rho^+$.

Now let
\begin{align*}
\tilde{\Omega}_\rho = (G/H) \times \mathbb{B}_\rho^+, \qquad \Omega_\rho = \Phi(\tilde{\Omega}_\rho),
\end{align*}
and decompose the boundary of $\Omega_\rho$ by setting
\begin{align*}
\partial_V \Omega_\rho &= \Phi\left((G/H) \times \partial_V \mathbb{B}_\rho^+\right), \\
\partial_L \Omega_\rho &= \Phi\left((G/H) \times \partial_L \mathbb{B}_\rho^+\right).
\end{align*}
The domain $\Omega_\rho$ need not be connected, since $G/H$ may be disconnected. Nevertheless, every component of $\Omega_\rho$ meets $\partial_V \Omega_\rho$, and $G$ permutes these components transitively, since it acts transitively on $G/H$. Thus Remark \ref{rem: disconnected problem} applies, and zero is a simple eigenvalue of the $G$-invariant Steklov--Neumann problem on $\Omega_\rho$.

Since $\Phi : \tilde{\Omega}_\rho \to \Omega_\rho$ is $G$-equivariant and preserves the boundary decomposition, Lemma \ref{quasi-isometry steklov-neumann} gives
\begin{align*}
C_\rho^{-(2n - 1)}\overline{\mu}_k^G(\tilde{\Omega}_\rho,g_{\mathrm{mod}}) \leq \overline{\mu}_k^G(\Omega_\rho,g_0) \leq C_\rho^{2n - 1}\overline{\mu}_k^G(\tilde{\Omega}_\rho,g_{\mathrm{mod}}).
\end{align*}
By Step 2,
\begin{align*}
\overline{\mu}_k^G(\tilde{\Omega}_\rho,g_{\mathrm{mod}}) = C_k \rho^{-\frac{2q}{n}},
\end{align*}
where $C_k > 0$ is independent of $\rho$, for every $k \geq 1$. Therefore,
\begin{align*}
C_k C_\rho^{-(2n - 1)}\rho^{-\frac{2q}{n}} \leq \overline{\mu}_k^G(\Omega_\rho,g_0) \leq C_k C_\rho^{2n - 1}\rho^{-\frac{2q}{n}}.
\end{align*}
Since $q > 0$ and $C_\rho \to 1$, we conclude that, for every $k \geq 1$,
\begin{align*}
\lim_{\rho \to 0}\overline{\mu}_k^G(\Omega_\rho,g_0) = +\infty.
\end{align*}

\textbf{Step 5: Conformal localisation}

\bigskip

It remains to transfer the divergence phenomenon obtained in Step 4 to the $G$-invariant Steklov spectrum of $M$. The idea is to keep the metric unchanged on $\Omega_\rho$ and collapse its complement conformally.

We first describe this construction in a more general setting. Fix a smooth $G$-invariant Riemannian metric $g$ on $M$, and let $\Omega \subset M$ be a compact $G$-invariant submanifold of codimension zero with Lipschitz boundary, not necessarily connected. Suppose that
\begin{align*}
\partial \Omega = \overline{\partial_V \Omega} \cup \overline{\partial_L \Omega},
\end{align*}
where $\partial_V \Omega$ and $\partial_L \Omega$ are disjoint relatively open smooth subsets of $\partial \Omega$, whose closures meet along a smooth codimension-two submanifold. We assume that every connected component of $\Omega$ meets $\partial_V \Omega$, that $G$ preserves both boundary portions, and that
\begin{align*}
\Omega \cap \partial M = \overline{\partial_V \Omega}, \qquad \partial \Omega \cap \interior M = \partial_L \Omega.
\end{align*}

For $0 < \varepsilon < 1$, define
\begin{align*}
f_\varepsilon =
\begin{cases}
1 & \text{on } \Omega,\\
\varepsilon^2 & \text{on } M \setminus \Omega,
\end{cases}
\qquad
g_\varepsilon = f_\varepsilon g.
\end{align*}
Although $g_\varepsilon$ is discontinuous across $\partial_L \Omega$, its Steklov eigenvalues can still be defined through the standard variational characterisation. More precisely, for $u \in H^1(M)^G$, the corresponding Rayleigh quotient is
\begin{align*}
R_\varepsilon(u) = \frac{ \displaystyle	\int_\Omega \vert \nabla u \vert_g^2 \, \dV_g +	\varepsilon^{n - 2}
\int_{M \setminus \Omega} \vert \nabla u \vert_g^2 \, \dV_g}
{ \displaystyle	\int_{\partial_V \Omega} u^2 \, \dA_g +	\varepsilon^{n - 1}	\int_{\partial M \setminus \overline{\partial_V \Omega}} u^2 \, \dA_g}.
\end{align*}
We denote by $\sigma_k^G(g_\varepsilon)$ the eigenvalues obtained by applying the usual min--max characterisation to this quotient. Since $n \geq 3$, as $\varepsilon \to 0$ the contribution of $M \setminus \Omega$ to the numerator vanishes, as does the contribution of $\partial M \setminus \overline{\partial_V \Omega}$ to the denominator. This suggests that the limiting problem is precisely the Steklov--Neumann problem on $\Omega$.

The next two propositions make this observation precise and allow us to replace the piecewise smooth metrics $g_\varepsilon$ by smooth ones. Their proofs are based on an argument of Enciso--Pistoia--Provenzano \cite{enciso-pistoia-provenzano}, which we adapt to the present $G$-invariant setting in Section \ref{sec: conformal localisation}.

\begin{prop}\label{spectral localisation}
Under the assumptions above, for every $k \geq 0$ for which the corresponding $G$-invariant Steklov--Neumann eigenvalue is defined, we have
\begin{align*}
\lim_{\varepsilon \to 0} \sigma_k^G(g_\varepsilon) = \mu_k^G(\Omega, g), \\
\lim_{\varepsilon \to 0} \overline{\sigma}_k^G(g_\varepsilon) = \overline{\mu}_k^G(\Omega, g).
\end{align*}
\end{prop}

We next approximate $g_\varepsilon$ by smooth conformal metrics.

\begin{prop} \label{prop: smooth approximation}
Under the assumptions above, fix $0 < \varepsilon < 1$. There exists a sequence of smooth $G$-invariant metrics $g_{\varepsilon, j}$ conformal to $g$ such that, for every $k \geq 0$ for which $\sigma_k^G(g_\varepsilon)$ is defined,
\begin{align*}
\lim_{j \to \infty} \sigma_k^G(g_{\varepsilon,j}) &= \sigma_k^G(g_\varepsilon), \\
\lim_{j \to \infty} \overline{\sigma}_k^G(g_{\varepsilon,j}) &= \overline{\sigma}_k^G(g_\varepsilon).
\end{align*}
Moreover, these metrics can be chosen in the form $g_{\varepsilon, j} = f_{\varepsilon, j} g$, where $\varepsilon^2 \leq f_{\varepsilon, j} \leq 1$ and $f_{\varepsilon, j} \to f_\varepsilon$ pointwise on $M$.
\end{prop}

We now apply the preceding results to the domains $\Omega_\rho$ constructed in Step 4, with $g = g_0$. By that step,
\begin{align*}
\lim_{\rho \to 0} \overline{\mu}_1^G(\Omega_\rho, g_0) = +\infty.
\end{align*}
Let $L > 0$. We may choose $\rho > 0$ sufficiently small so that
\begin{align*}
\overline{\mu}_1^G(\Omega_\rho,g_0) > 2L.
\end{align*}
For this fixed $\rho$, let $g_{\rho,\varepsilon}$ denote the piecewise conformal metric obtained from $\Omega_\rho$ as above. Proposition \ref{spectral localisation} and the volume identities imply that
\begin{align*}
\lim_{\varepsilon \to 0} \overline{\sigma}_1^G(g_{\rho, \varepsilon}) = \overline{\mu}_1^G(\Omega_\rho, g_0).
\end{align*}
We may therefore choose $\varepsilon > 0$ sufficiently small so that
\begin{align*}
\overline{\sigma}_1^G(g_{\rho, \varepsilon}) > \frac{3L}{2}.
\end{align*}
Finally, Proposition \ref{prop: smooth approximation} gives a sequence of smooth $G$-invariant metrics $g_{\rho, \varepsilon, j}$ conformal to $g_0$ such that
\begin{align*}
\lim_{j \to \infty} \overline{\sigma}_1^G(g_{\rho, \varepsilon,j}) = \overline{\sigma}_1^G(g_{\rho, \varepsilon}).
\end{align*}
Hence, for $j$ sufficiently large,
\begin{align*}
\overline{\sigma}_1^G(g_{\rho,\varepsilon,j}) > L.
\end{align*}
Since $L > 0$ is arbitrary and every $g_{\rho,\varepsilon, j}$ is a smooth $G$-invariant metric conformal to $g_0$, we conclude that
\begin{align*}
\sup_g \sigma_1^G(g) \Vol(M,g)^{\frac{2 - n}{n}} \Vol(\partial M, g) =	+\infty.
\end{align*}
This proves Theorem \ref{thm: infinite sup}.

\medskip

We conclude this section by proving Corollary \ref{cor: first eigenfunctions}.

\begin{proof}[Proof of Corollary~\ref{cor: first eigenfunctions}]
Fix $K \geq 1$. By Hassannezhad's upper bound \cite[Theorem~4.1]{hassannezhad} mentioned in the Introduction, there exists $C > 0$ such that $\overline{\sigma}_K(g) \leq C K^{2/n}$ for every $g \in [g_0]$. By Theorem \ref{thm: infinite sup}, we may choose a smooth $G$-invariant metric $g_K$ conformal to $g_0$ such that 
\begin{align*}
\overline{\sigma}_K(g_K) \leq CK^{2/n} <\overline{\sigma}_1^G(g_K).
\end{align*}
Since the normalisation factors coincide, this implies $\sigma_K(g_K) < \sigma_1^G(g_K)$. Consequently, no nonzero	eigenfunction associated with any of the eigenvalues $\sigma_1(g_K), \dots, \sigma_K(g_K)$ is $G$-invariant.
\end{proof}

\section{Conformal localisation} \label{sec: conformal localisation}

In this section, we establish the spectral convergence and smooth approximation results used in the conformal localisation argument of the preceding section.

\subsection{Spectral convergence under conformal collapse} \label{sec: Spectral convergence under conformal collapse}

We now prove Proposition \ref{spectral localisation}. The argument is based on the ideas of Enciso--Pistoia--Provenzano \cite{enciso-pistoia-provenzano}. Since we only need convergence of the eigenvalues, the proof can be formulated directly in terms of the min--max characterisations.

Throughout this section, $C > 0$ denotes a constant whose value may change from line to line and which is independent of the parameters tending to zero, unless otherwise specified.

We shall repeatedly use the following Poincar\'e--trace inequality. Let $(Y, g)$ be a compact Lipschitz Riemannian manifold with boundary, not necessarily connected, and let $\Sigma \subset \partial Y$ be relatively open. Suppose that $\Sigma$ meets every connected component of $Y$. Then there exists $C > 0$ such that
\begin{align*}
\norm{u}_{H^1(Y)}^2 \leq C \left( \int_Y \vert \nabla u \vert^2 \, \dV_g + \int_\Sigma u^2 \, \dA_g \right)
\end{align*}
for every $u \in H^1(Y)$.

For connected $Y$, this is the equivalence of norms recalled in \cite[Appendix B]{enciso-pistoia-provenzano}. The general case follows by applying this inequality to each connected component, using the corresponding portion of $\Sigma$, and summing over the finitely many components.

\begin{proof}[Proof of Proposition \ref{spectral localisation}]
Set
\begin{align*}
D = \overline{M \setminus \Omega}, \qquad \Gamma = \partial M \setminus \overline{\partial_V \Omega}.
\end{align*}
For $u, v \in H^1(M)^G$, define
\begin{align*}
\mathcal{E}_\varepsilon(u, v) &= \int_\Omega g( \nabla u, \nabla v) \, \dV_g + \varepsilon^{n - 2} \int_D g( \nabla u, \nabla v) \, \dV_g, \\
\mathcal{B}_\varepsilon(u, v) &= \int_{\partial_V \Omega} u v \, \dA_g + \varepsilon^{n - 1} \int_\Gamma u v \, \dA_g.
\end{align*}
For brevity, write $\mathcal{E}_\varepsilon(u) = \mathcal{E}_\varepsilon(u, u)$ and $\mathcal{B}_\varepsilon(u) = \mathcal{B}_\varepsilon(u, u)$. Thus, the Rayleigh quotient associated with $g_\varepsilon$ is
\begin{align*}
R_\varepsilon(u) = \frac{ \mathcal E_\varepsilon(u)} { \mathcal B_\varepsilon(u) }.
\end{align*}

We will prove that, for every fixed $k \geq 0$, there exists $C > 0$, independent of $\varepsilon$, such that
\begin{align*}
(1 - C \varepsilon) \mu_k^G(\Omega, g) \leq \sigma_k^G(g_\varepsilon) \leq \mu_k^G(\Omega, g) + C \varepsilon^{n - 2}
\end{align*}
for all sufficiently small $\varepsilon > 0$.
	
\medskip
\noindent \textbf{Upper bound.} Let $w_0, \ldots, w_k \in H^1(\Omega)^G$ be Steklov--Neumann eigenfunctions corresponding to $\mu_0^G(\Omega, g), \ldots, \mu_k^G(\Omega, g)$, and assume their traces are orthonormal in $L^2(\partial_V \Omega, g)$. Let $W_k$ be the span of $w_0, \dots, w_k$ in $H^1(\Omega)^G$. Then, for every $w \in W_k$,
\begin{align*}
\int_\Omega \vert \nabla w \vert^2 \, \dV_g \leq \mu_k^G(\Omega, g) \int_{\partial_V \Omega} w^2 \, \dA_g.
\end{align*}
By the Poincar\'e--trace inequality on $\Omega$, we also have
\begin{align*}
\norm{w}_{H^1(\Omega)}^2 \leq C \int_{\partial_V \Omega} w^2 \, \dA_g.
\end{align*}
	
Since $\Omega$ has Lipschitz boundary, there exists a bounded extension operator $E_0 : H^1(\Omega) \to H^1(M)$. 
As $E_0$ need not preserve $G$-invariance, for $w \in H^1(\Omega)^G$ we define
\begin{align*}
Ew = \int_G \gamma \cdot E_0 w \, \dd \mu(\gamma),
\end{align*}
where $\mu$ is the normalised Haar measure on $G$. Since $\Omega$ is $G$-invariant, $Ew$ is a $G$-invariant extension of $w$, so $E : H^1(\Omega)^G \to H^1(M)^G$ is a well-defined extension operator. It is also bounded. Indeed, the norm on $H^1(M)$ is a convex function and $\mu$ is a probability measure, so Jensen's inequality gives
\begin{align*}
\norm{Ew}_{H^1(M)} &\leq \int_G \norm{\gamma \cdot E_0 w}_{H^1(M)} \, \dd \mu(\gamma) = \int_G \norm{E_0 w}_{H^1(M)} \, \dd \mu(\gamma) = \norm{E_0 w}_{H^1(M)} \leq C \norm{w}_{H^1(\Omega)},
\end{align*}
where we used the fact that the $G$-action is isometric and that $E_0$ is bounded.

Consequently, for every $w \in W_k$,
\begin{align*}
\mathcal E_\varepsilon(Ew) &\leq \int_\Omega \vert \nabla w \vert^2 \, \dV_g + C \varepsilon^{n - 2} \norm{w}_{H^1(\Omega)}^2 \leq \left( \mu_k^G(\Omega, g) + C \varepsilon^{n - 2} \right) \int_{\partial_V \Omega} w^2 \, \dA_g,
\end{align*}
and
\begin{align*}
\mathcal B_\varepsilon(Ew) &\geq \int_{\partial_V \Omega} w^2 \, \dA_g.
\end{align*}
Since $E$ is injective on $W_k$, the space $E(W_k)$ has dimension $k + 1$. The min--max characterisation then gives
\begin{align}
\sigma_k^G(g_\varepsilon) \leq \mu_k^G(\Omega, g) + C \varepsilon^{n - 2}. \label{eq:conformal-collapse-upper}
\end{align}
	
\medskip
\noindent \textbf{Lower bound.} Let $u_{0, \varepsilon}, \dots, u_{k, \varepsilon}$ be $G$-invariant eigenfunctions associated with $\sigma_0^G(g_\varepsilon), \dots, \sigma_k^G(g_\varepsilon)$, and assume they are orthonormal with respect to the bilinear form $\mathcal{B}_\varepsilon$.

Set $F_{\varepsilon,k} = \Span \{u_{0,\varepsilon}, \ldots, u_{k,\varepsilon}\}$, so that $\dim F_{\varepsilon, k} = k + 1$. Moreover, for every $u = \sum_{i=0}^k a_i u_{i,\varepsilon} \in F_{\varepsilon,k}$, the weak eigenvalue equations and orthonormality give
\begin{align*}
\mathcal{B}_\varepsilon(u) = \sum_{i=0}^k a_i^2, \qquad \mathcal{E}_\varepsilon(u) = \sum_{i=0}^k \sigma_i^G(g_\varepsilon) a_i^2 \leq \sigma_k^G(g_\varepsilon) \mathcal{B}_\varepsilon(u).
\end{align*}
By \eqref{eq:conformal-collapse-upper}, $\sigma_k^G(g_\varepsilon)$ is uniformly bounded for sufficiently small $\varepsilon$. Thus,
\begin{align}
\mathcal{E}_\varepsilon(u) \leq \sigma_k^G(g_\varepsilon) \mathcal{B}_\varepsilon(u) \leq C \mathcal{B}_\varepsilon(u) \label{eq:conformal-collapse-energy}
\end{align}
for every $u \in F_{\varepsilon,k}$ and all sufficiently small $0 <\varepsilon < 1$.

Applying the Poincar\'e--trace inequality recalled above to $Y = M$ and $\Sigma = \partial_V \Omega$, and using the continuity of the trace map, we obtain
\begin{align*}
\int_\Gamma u^2 \, \dA_g \leq C \left( \int_M \vert \nabla u \vert^2 \, \dV_g + \int_{\partial_V \Omega} u^2 \, \dA_g \right).
\end{align*}
Since $0 < \varepsilon < 1$ and $n \geq 3$, we have $\int_M \vert \nabla u \vert^2 \, \dV_g \leq \varepsilon^{2 - n} \mathcal E_\varepsilon(u)$. Hence, for every $u \in F_{\varepsilon,k}$, \eqref{eq:conformal-collapse-energy} implies that
\begin{align*}
\varepsilon^{n - 1} \int_\Gamma u^2 \, \dA_g &\leq C \varepsilon \mathcal E_\varepsilon(u) + C \varepsilon^{n - 1} \mathcal B_\varepsilon(u) \\
&\leq C \varepsilon \mathcal B_\varepsilon(u).
\end{align*}
It follows that
\begin{align}
\int_{\partial_V \Omega} u^2 \, \dA_g \geq (1 - C \varepsilon) \mathcal B_\varepsilon(u) \label{eq:conformal-collapse-boundary}
\end{align}
for every $u \in F_{\varepsilon,k}$.
	
For $\varepsilon$ sufficiently small, \eqref{eq:conformal-collapse-boundary} implies that restriction to $\Omega$ is injective on $F_{\varepsilon,k}$. Indeed, if $u \in F_{\varepsilon,k}$ vanishes on $\Omega$, then \eqref{eq:conformal-collapse-boundary} implies $\mathcal{B}_\varepsilon(u) = 0$, hence $u = 0$ by the orthonormality of the chosen eigenfunctions. Thus,
\begin{align*}
V_{\varepsilon,k} = \{ u \vert_\Omega : u \in F_{\varepsilon,k} \} \subset H^1(\Omega)^G
\end{align*}
has dimension $k + 1$. By the min--max characterisation of the $G$-invariant Steklov--Neumann eigenvalues, \eqref{eq:conformal-collapse-energy}, and \eqref{eq:conformal-collapse-boundary},
\begin{align*}
\mu_k^G(\Omega, g) &\leq \max_{0 \neq u \in F_{\varepsilon,k}} \frac{\displaystyle \int_\Omega \vert \nabla u \vert^2 \, \dV_g}{\displaystyle \int_{\partial_V \Omega} u^2 \, \dA_g} \\
&\leq \frac{\sigma_k^G(g_\varepsilon)}{1 - C \varepsilon}.
\end{align*}
	
Combining these bounds yields the desired convergence of the eigenvalues as $\varepsilon \to 0$. Finally, since $g_\varepsilon = g$ on $\Omega$ and $g_\varepsilon = \varepsilon^2 g$ on $M \setminus \Omega$, we have
\begin{align*}
\Vol(M, g_\varepsilon) &= \Vol(\Omega, g) + \varepsilon^n \Vol(D, g) \to \Vol(\Omega, g), \\
\Vol(\partial M, g_\varepsilon) &= \Vol(\partial_V \Omega, g) + \varepsilon^{n - 1} \Vol(\Gamma, g) \to \Vol(\partial_V \Omega, g).
\end{align*}
This implies
\begin{align*}
\lim_{\varepsilon \to 0} \overline{\sigma}_k^G(g_\varepsilon) = \overline{\mu}_k^G(\Omega, g),
\end{align*}
as desired.
\end{proof}

\subsection{Smooth approximation of the collapsing metrics}

We now prove Proposition \ref{prop: smooth approximation}, following the argument of \cite[Proposition 2.4]{enciso-pistoia-provenzano}.
\begin{proof}[Proof of Proposition \ref{prop: smooth approximation}]
Fix $0 < \varepsilon < 1$. For each $j \geq 1$, choose a smooth $G$-invariant function $\chi_j : M \to [0,1]$, equal to $1$ on $\Omega$ and supported in $\Omega_j = \{ x \in M : d_g(x, \Omega) < 1/j \}$. Such functions exist by averaging over $G$, since $\Omega$ and $\Omega_j$ are $G$-invariant. Set
\begin{align*}
f_{\varepsilon,j} = \varepsilon^2 + (1 - \varepsilon^2) \chi_j, \qquad g_{\varepsilon, j} = f_{\varepsilon, j}g.
\end{align*}
Then $g_{\varepsilon,j}$ is smooth and $G$-invariant, $f_\varepsilon \leq f_{\varepsilon, j} \leq 1$, and $f_{\varepsilon, j} \to f_\varepsilon$ pointwise on $M$.

Exactly as in the proof of \cite[Proposition 2.4]{enciso-pistoia-provenzano}, the convergence of the eigenvalues follows from Colin de Verdi\`ere's approximation theorem \cite[Theorem I.8]{colin_verdiere}. Indeed, the same argument applies to the $G$-invariant subspace $H^{1/2}(\partial M)^G$, since all the metrics are $G$-invariant and harmonic extension preserves $G$-invariance. Thus, for every $k \geq 0$,
\begin{align*}
\sigma_k^G(g_{\varepsilon,j}) \longrightarrow \sigma_k^G(g_\varepsilon) \qquad \text{as } j \to \infty.
\end{align*}
	
Finally, dominated convergence gives
\begin{align*}
\Vol(M, g_{\varepsilon,j}) \to \Vol(M, g_\varepsilon), \qquad \Vol(\partial M, g_{\varepsilon,j}) \to \Vol(\partial M, g_\varepsilon), \qquad \text{as } j \to \infty.
\end{align*}
Together with the eigenvalue convergence, this yields $\overline{\sigma}_k^G(g_{\varepsilon,j}) \to \overline{\sigma}_k^G(g_\varepsilon)$, as desired.
\end{proof}

\section{Euclidean domains of revolution} \label{sec: Revolution-domains}

In this section, we prove Theorem \ref{thm: revolution balls unbounded}. We give the details in dimension three, where $G = \mathrm{O}(2)$ acts on the first two coordinates by orthogonal transformations, and explain the extension to higher dimensions at the end of the proof. The construction starts from cylinders of increasing radius and fixed height. We round their edges and introduce increasingly rapid oscillations on a central portion of the lateral boundary.

We begin with a Poincar\'e-type inequality for functions with zero mean on a portion of the boundary.

\begin{lem} \label{lem:rev:boundary-mean-poincare}
Let $(Y, g)$ be a compact connected Lipschitz Riemannian manifold with boundary, and let $\Sigma \subset \partial Y$ be nonempty and relatively open. Then there exists $C > 0$ such that
\begin{align*}
\norm{u}_{H^1(Y)}^2 \leq C \int_Y \vert \nabla u \vert^2 \, \dV_g
\end{align*}
for every $u \in H^1(Y)$ satisfying $\int_\Sigma u \, \dA_g = 0$.
\end{lem}

\begin{proof}
Suppose the conclusion is false. Then there exists a sequence $u_j \in H^1(Y)$ such that
\begin{align*}
\norm{u_j}_{H^1(Y)} = 1, \qquad \int_\Sigma u_j \dA_g = 0, \qquad \int_Y \vert \nabla u_j \vert^2 \dV_g \to 0.
\end{align*}
By weak compactness in $H^1(Y)$ and compactness of the embedding $H^1(Y) \hookrightarrow L^2(Y)$, after passing to a subsequence, we may assume that $u_j \rightharpoonup u$ weakly in $H^1(Y)$ and $u_j \to u$ strongly in $L^2(Y)$, for some $u \in H^1(Y)$. Since $\nabla u_j \to 0$ in $L^2(Y)$, we have $\nabla u = 0$, and hence $u$ is constant because $Y$ is connected. Moreover, by continuity of the trace operator, the zero-mean condition passes to the limit, giving $\int_\Sigma u \dA_g = 0$. Since $\Sigma$ has positive measure in $\partial Y$, we conclude that $u = 0$. Thus $u_j \to 0$ strongly in $H^1(Y)$, contradicting $\norm{u_j}_{H^1(Y)} = 1$.
\end{proof}

We next use this lemma to obtain two estimates for invariant functions on cylinders. For $L > 0$, denote by $\mathbb{B}_L^2 = \{ x \in \mathbb R^2 : \vert x \vert < L\}$ the Euclidean disc of radius $L$, and set $Q_L = \overline{\mathbb{B}_L^2} \times [-1, 1]$.

\begin{lem} \label{lem: poincare on cylinders}
There exists $C > 0$, independent of $L \geq 1$, such that the following holds. Let $w \in H^1(Q_L)^G$ and write $w(r, \theta, z) = v(r, z)$ in cylindrical coordinates. If its lateral trace satisfies $\int_{-1}^1 v(L, z) \, \dd z = 0$, then
\begin{align}
L \int_{-1}^1 v(L, z)^2 \, \dd z &\leq C \int_{Q_L} \vert \nabla w \vert^2 \, \dd V, \label{eq:rev:cylinder-lateral} \\
\int_{\overline{\mathbb{B}_L^2} \times \{-1, 1\}} w^2 \, \dd A &\leq C L^2 \int_{Q_L} \vert \nabla w \vert^2 \, \dd V. \label{eq:rev:cylinder-faces}
\end{align}
\end{lem}

\begin{proof}
We first assume that $w$ is smooth and $G$-invariant, and prove the lateral trace estimate. Consider the rectangle $S_L = [L - 1/2, L] \times [-1, 1]$ and its right-hand side $\Sigma_L = \{L\} \times (-1, 1)$. Since the pairs $(S_L, \Sigma_L)$ are all isometric by translations, the continuity of the trace map $H^1(S_L) \to L^2(\Sigma_L)$ and Lemma \ref{lem:rev:boundary-mean-poincare} give
\begin{align*}
\int_{-1}^1 v(L, z)^2 \, \dd z \leq C \norm{v}_{H^1(S_L)}^2 \leq C \int_{S_L} (v_r^2 + v_z^2) \, \dd r \, \dd z,
\end{align*}
where $C$ is independent of $L$. Moreover, $r \geq L/2$ on $S_L$ whenever $L \geq 1$. So, multiplying this inequality by $L$ and using the cylindrical expression for the Dirichlet energy, we obtain
\begin{align*}
L \int_{-1}^1 v(L, z)^2 \, \dd z &\leq C \int_{L - 1/2}^L \int_{-1}^1 (v_r^2 + v_z^2) r \, \dd z \, \dd r \\
&\leq C \int_{Q_L} \vert \nabla w \vert^2 \, \dd V,
\end{align*}
which proves \eqref{eq:rev:cylinder-lateral}.
	
We now estimate the $L^2$-norm in the cylinder. For $0 < r < L$, the fundamental theorem of calculus gives
\begin{align*}
v(r, z) = v(L, z) - \int_r^L v_s(s, z) \, \dd s.
\end{align*}
By Cauchy--Schwarz, we have:
\begin{align*}
\left\vert \int_r^L v_s(s, z) \, \dd s \right\vert^2 &\leq \left( \int_r^L \vert v_s(s, z) \vert^2 s \, \dd s \right) \left( \int_r^L \frac{1}{s} \, \dd s \right) \\
&\leq \log \frac{L}{r} \int_0^L \vert v_s(s, z) \vert^2 s \, \dd s.
\end{align*}
So, using $(a + b)^2 \leq 2a^2 + 2b^2$, we obtain
\begin{align*}
\vert v(r, z) \vert^2 \leq 2 v(L, z)^2 + 2 \log \frac{L}{r} \int_0^L \vert v_s(s, z) \vert^2 s \, \dd s.
\end{align*}
Thus, multiplying this inequality by $2 \pi r$ and integrating over $(0, L) \times (-1, 1)$, we find
\begin{align*}
\int_{Q_L} w^2 \, \dd V &\leq 4 \pi \left( \int_0^L r \, \dd r \right) \int_{-1}^1 v(L, z)^2 \, \dd z \\
&\quad + 4 \pi \left( \int_0^L r \log \frac{L}{r} \, \dd r \right) \int_{-1}^1 \int_0^L \vert v_s(s, z) \vert^2 s \, \dd s \, \dd z \\
&= 2 \pi L^2 \int_{-1}^1 v(L, z)^2 \, \dd z + \frac{L^2}{2} \int_{Q_L} \vert \partial_r w \vert^2 \, \dd V.
\end{align*}
Now, the lateral trace estimate \eqref{eq:rev:cylinder-lateral} implies
\begin{align*}
2 \pi L^2 \int_{-1}^1 v(L, z)^2 \, \dd z \leq C L \int_{Q_L} \vert \nabla w \vert^2 \, \dd V \leq C L^2 \int_{Q_L} \vert \nabla w \vert^2 \, \dd V,
\end{align*}
where we used $L \geq 1$. Since $\vert \partial_r w \vert^2 \leq \vert \nabla w \vert^2$, we conclude that
\begin{align*} 
\int_{Q_L} w^2 \, \dd V \leq C L^2 \int_{Q_L} \vert \nabla w \vert^2 \, \dd V.
\end{align*}

Finally, the standard one-dimensional trace inequality on $(-1, 1)$ gives, for every $r \in (0, L)$,
\begin{align*}
v(r, 1)^2 + v(r, -1)^2 \leq C \int_{-1}^1 (v(r, z)^2 + v_z(r, z)^2) \, \dd z.
\end{align*}
Multiplying by $2 \pi r$, integrating and using the preceding interior $L^2$-estimate yields
\begin{align*}
\int_{\overline{\mathbb{B}_L^2} \times \{-1, 1\}} w^2 \, \dd A \leq C \left( \int_{Q_L} w^2 \, \dd V + \int_{Q_L} \vert \partial_z w \vert^2 \, \dd V \right) \leq C L^2 \int_{Q_L} \vert \nabla w \vert^2 \, \dd V.
\end{align*}
This proves \eqref{eq:rev:cylinder-faces}. 
	
We finally remove the smoothness assumption. Let $w \in H^1(Q_L)^G$ satisfy the zero-mean condition on the lateral boundary. By density, choose smooth functions $\varphi_j$ converging to $w$ in $H^1(Q_L)$. Averaging over $G$, we obtain smooth $G$-invariant functions $w_j$ converging to $w$ in $H^1(Q_L)$. Writing $w_j(r, \theta, z) = v_j(r, z)$, set
\begin{align*}
\alpha_j = \frac{1}{2} \int_{-1}^1 v_j(L, z) \, \dd z, \qquad \tilde w_j = w_j - \alpha_j.
\end{align*}
By continuity of the trace map and the zero-mean condition on $w$, we have $\alpha_j \to 0$. Hence $\tilde w_j \to w$ in $H^1(Q_L)$, and each $\tilde w_j$ is smooth, $G$-invariant, and has zero lateral mean. Applying the preceding estimates to $\tilde w_j$ and passing to the limit proves the result.	
\end{proof}

\begin{proof}[Proof of Theorem \ref{thm: revolution balls unbounded}]

We construct compact smooth $G$-invariant domains $\Omega_R \subset \mathbb R^3$ for integers $R \geq 2$. Figures~\ref{profile} and~\ref{OmegaR} illustrate the construction for $R = 2$, showing, respectively, the meridian profile in the $rz$-plane and the solid domain obtained by revolution about the $z$-axis.

Start with the solid cylinder $\overline{\mathbb B_R^2} \times [-1, 1]$. In cylindrical coordinates $(r, \theta, z)$, its meridian region is the rectangle $[0, R] \times [-1, 1]$ in the half-plane $r \geq 0$. We round the two right-hand corners inward and replace a central portion of the right-hand vertical side by an oscillating curve. The rounding profile will be the same for every $R$.

\begin{figure}[!htbp]
\centering
\includegraphics[width=0.75\linewidth]{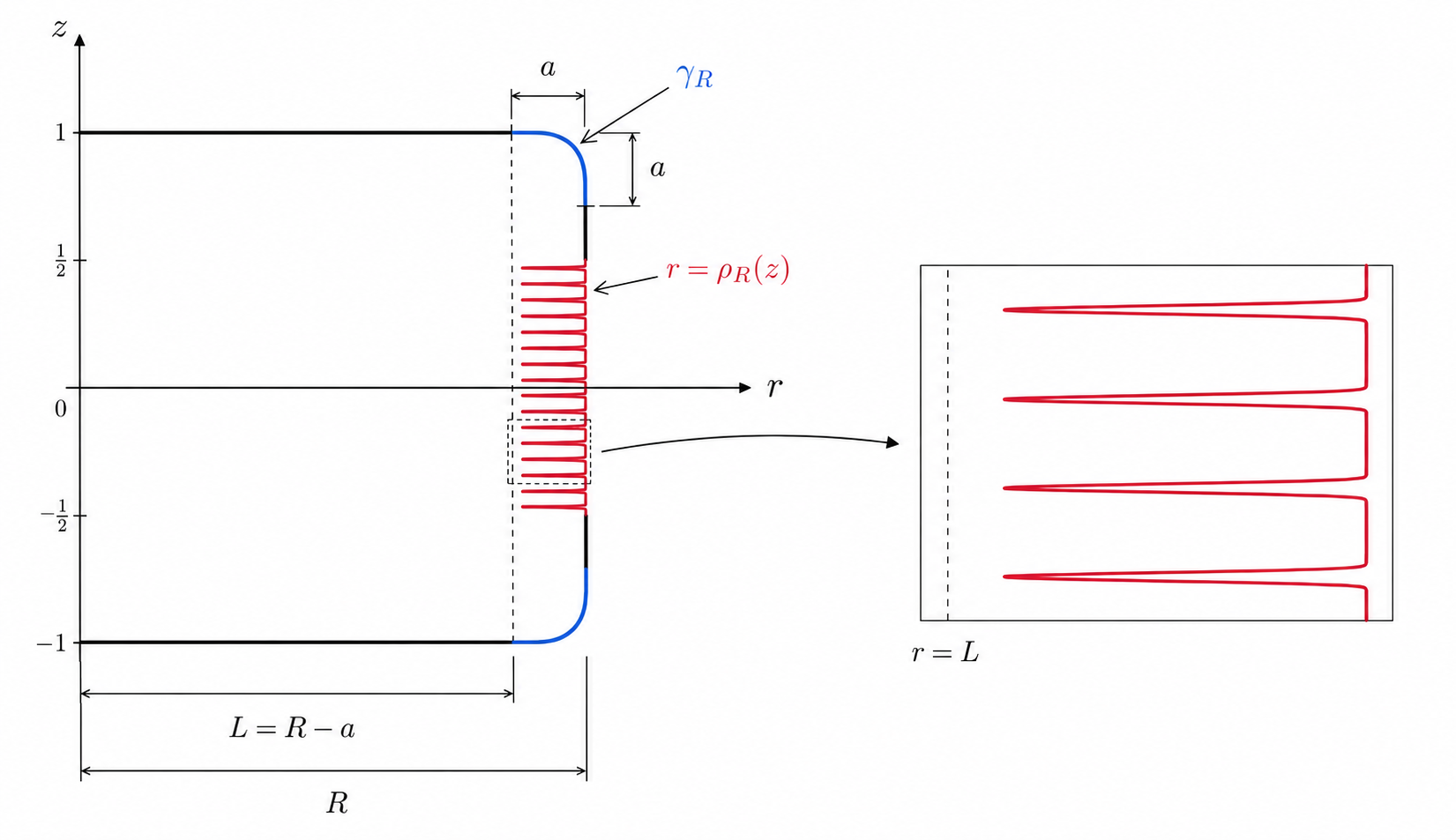}
\caption{The meridian profile of $\Omega_R$ for $R = 2$. The blue curves are the rounded arcs and the red curve is $r = \rho_R(z)$. The dashed line $r = L$ identifies the inner cylinder. The detail on the right enlarges four consecutive oscillations of the profile curve.}
\label{profile}
\end{figure}
	
To define the rounded corners, let
\begin{align*}
p(t) = 
\begin{cases}
e^{-1/t}, & t > 0, \\
0, & t \leq 0, 
\end{cases}
\qquad \chi(t) = \frac{p(t)}{p(t) + p(1 - t)}, \qquad 0 \leq t \leq 1.
\end{align*}
The function $\chi$ is smooth and strictly increasing, with $\chi(0) = 0$, $\chi(1) = 1$ and $\chi(1 - t) = 1 - \chi(t)$. Moreover, all its derivatives of positive order vanish at both endpoints. Set
\begin{align*}
\xi(t) = \frac{\pi}{2} \chi(t), \qquad a = \frac{1}{2} \int_0^1 \cos \xi(t) \, \dd t = \frac{1}{2} \int_0^1 \sin \xi(t) \, \dd t.
\end{align*}
The equality follows from the symmetry of $\chi$. Since $0 < \xi(t) < \pi/2$ for $0 < t < 1$, we have $1 < \cos \xi(t) + \sin \xi(t)$ and $\cos \xi(t) < 1$ there. Consequently, $1/4 < a < 1/2$.
	
For each integer $R \geq 2$, put $L = R - a$ and define the upper rounded arc by
\begin{align*}
\gamma_R(t) = (L, 1) + \frac{1}{2} \int_0^t (\cos \xi(s), - \sin \xi(s)) \, \dd s, \qquad 0 \leq t \leq 1.
\end{align*}
Since $\gamma_R'(t) = \frac{1}{2} (\cos \xi(t), - \sin \xi(t))$, its radial coordinate increases, its height decreases, and its endpoints are $\gamma_R(0) = (L, 1)$ and $\gamma_R(1) = (R, 1 - a)$. The flatness of $\chi$ at the endpoints makes the arc join smoothly to the horizontal segment $z = 1$ and to the vertical segment $r = R$. Reflecting $\gamma_R$ across $z = 0$ gives the lower rounded arc.
	
For the oscillating portion, define
\begin{align*}
\eta_0(t) = \begin{cases} \displaystyle \exp \left( 16 - \frac{1}{(t - 1/4)(3/4 - t)} \right), & 1/4 < t < 3/4, \\ 0, & 0 \leq t \leq 1/4 \text{ or } 3/4 \leq t \leq 1, \end{cases}
\end{align*}
and let $\eta$ be its $1$-periodic extension. Thus $\eta \in C^\infty(\mathbb R)$ vanishes near every integer and satisfies $0 \leq \eta \leq 1$, with $\eta(1/2) = 1$. Let $\rho_R$ be the smooth function given by
\begin{align*}
\rho_R(z) = R - \frac{1}{R^2} \eta \left( R^4 \left( z + \frac{1}{2} \right) \right), \qquad \vert z \vert \leq 1/2
\end{align*}
Since $R^4$ is an integer, the graph of $\rho_R$ contains exactly $R^4$ periods and agrees with $r = R$ near $z = \pm 1/2$. Because $R^{-2} \leq 1/4 < a$, we have $L = R - a < R - R^{-2} \leq \rho_R(z) \leq R$, so the oscillations remain strictly to the right of $r = L$. Moreover, $a < 1/2$ ensures that the unchanged vertical segments $r = R$, $1/2 \leq \vert z \vert \leq 1 - a$, have positive length.
	
Join the two horizontal segments $z = \pm 1$, $0 \leq r \leq L$, the two rounded arcs, the unchanged vertical segments and the graph of $\rho_R$. Together with the rotation axis, these curves bound a compact meridian region in the $rz$ plane. Rotating this region around the $z$-axis defines the compact smooth domain $\Omega_R \subset \mathbb R^3$. The boundary of $\Omega_R$ is smooth by construction, since the curves join smoothly at their endpoints; where it meets the axis, it is locally a flat disc.
	
The domain $\Omega_R$ is diffeomorphic to a closed ball. Indeed, let $\Omega_R^0$ denote the closed rounded cylinder before introducing the oscillations. Define $F_R : \mathbb R^3 \to \mathbb R^3$ by
\begin{align*}
F_R(x, z) = (\alpha_R(z) x, z), \qquad \alpha_R(z) = \begin{cases} \rho_R(z)/R, & \vert z \vert \leq 1/2, \\ 1, & \vert z \vert > 1/2, \end{cases}
\end{align*}
where $(x, z) \in \mathbb R^2 \times \mathbb R$. Since $\alpha_R$ is smooth and strictly positive, $F_R$ is a smooth diffeomorphism satisfying $F_R(\Omega_R^0) = \Omega_R$.
	
\begin{figure}[!htbp]
\centering
\includegraphics[width=0.4\linewidth]{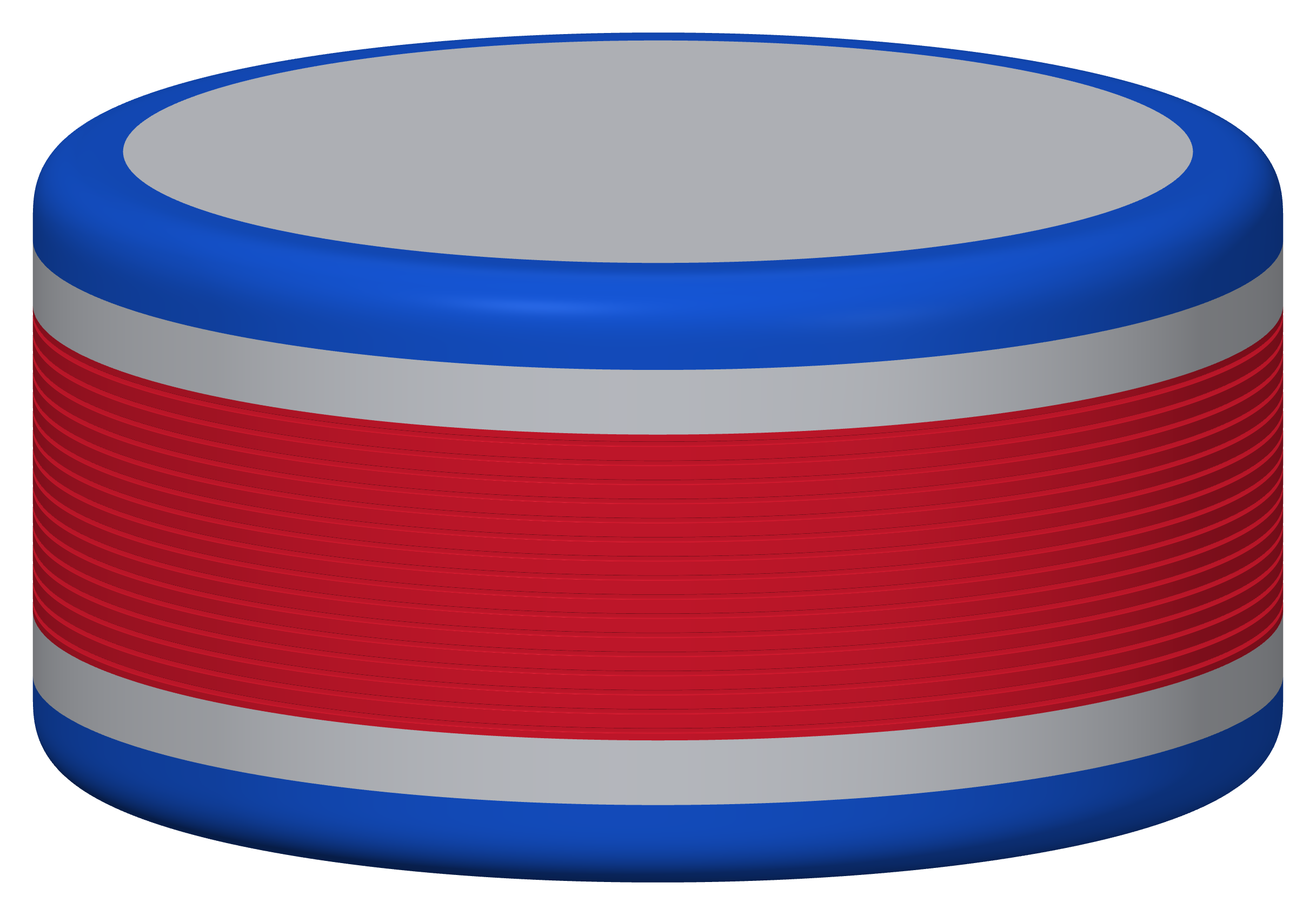}
\caption{The domain $\Omega_2$ obtained by revolving the meridian profile in Figure~\ref{profile} around the $z$-axis. The blue rounded parts, red oscillating band and unmodified cylinder portions correspond to those in the meridian construction. The top and bottom faces are flat discs.}
\label{OmegaR}
\end{figure}

For each $z \in [-1, 1]$, let $\lambda_R(z)$ be the outer radius of the meridian region at height $z$. In particular, $\lambda_R = \rho_R$ on $\vert z \vert \leq 1/2$, and $\lambda_R = R$ on the unchanged vertical portions. The monotonicity of the rounded arcs and the preceding bounds imply that $L \leq \lambda_R(z) \leq R$. Thus
\begin{align}
Q_L \subset \Omega_R \subset \overline {\mathbb{B}_R^2} \times [-1, 1]. \label{eq:rev:cylinder-inclusion}
\end{align}

We now establish an upper bound for the volume of $\Omega_R$, a lower bound for its boundary area, and a lower bound for $\sigma_1^G(\Omega_R)$. As usual, throughout the proof, $c, C > 0$ denote constants independent of $R$ and of the functions involved, whose values may change from line to line.
	
\medskip
\noindent \textbf{Volume upper bound.}
By \eqref{eq:rev:cylinder-inclusion}, we immediately obtain
\begin{align}
\Vol(\Omega_R) \leq \Vol( \overline{ \mathbb{B}_R^2} \times [-1, 1]) = 2 \pi R^2. \label{eq:rev:volume-bound}
\end{align}	
	
\medskip
\noindent \textbf{Boundary area lower bound.}
Let $\Gamma_R^{\mathrm{osc}}$ be the surface generated by the oscillating graph. It is parametrised by
\begin{align*}
X_R(\theta, z) = (\rho_R(z) \cos \theta, \rho_R(z) \sin \theta, z), \qquad 0 \leq \theta < 2 \pi, \quad -1/2 \leq z \leq 1/2.
\end{align*}
Consequently,
\begin{align*}
\Vol(\Gamma_R^{\mathrm{osc}}) = 2 \pi \int_{-1/2}^{1/2} \rho_R(z) \sqrt{1 + \rho_R'(z)^2} \, \dd z.
\end{align*}
Since 
\begin{align*}
\rho_R(z) \geq R - \frac{1}{R^2} \geq \frac{R}{2}, \qquad \rho_R'(z) = - R^2 \eta' \left( R^4 \left( z + \frac{1}{2} \right) \right),
\end{align*}
we obtain
\begin{align*}
\Vol(\Gamma_R^{\mathrm{osc}}) \geq \pi R \int_{-1/2}^{1/2} \vert \rho_R'(z) \vert \, \dd z = \pi R^3 \int_{-1/2}^{1/2} \left\vert \eta' \left( R^4 \left( z + \frac{1}{2} \right) \right) \right\vert \, \dd z.
\end{align*}
As $\eta$ is $1$-periodic and $R^4$ is an integer, the change of variables $t = R^4(z + 1/2)$ gives
\begin{align*}
\int_{-1/2}^{1/2} \left\vert \eta' \left( R^4 \left( z + \frac{1}{2} \right) \right) \right\vert \, \dd z = \int_0^1 \vert \eta'(t) \vert \, \dd t
\end{align*}
So, writing $I_\eta = \int_0^1 \vert \eta'(t) \vert \, \dd t > 0$, we conclude that
\begin{align}
\Vol(\partial \Omega_R) \geq \Vol(\Gamma_R^{\mathrm{osc}}) \geq \pi I_\eta R^3 = c R^3. \label{eq:rev:area-bound}
\end{align}
	
\medskip

\noindent \textbf{Lower bound for $\sigma_1^G(\Omega_R)$.}
We prove that every $u \in H^1(\Omega_R)^G$ with zero boundary mean satisfies
\begin{align}
\int_{\partial \Omega_R} u^2 \, \dd A \leq C R^2 \int_{\Omega_R} \vert \nabla u \vert^2 \, \dd V. \label{eq:rev:target}
\end{align}
This gives the required lower bound by the variational characterisation of $\sigma_1^G(\Omega_R)$. By smooth approximation, followed by averaging over $G$ and subtracting the boundary mean, it is enough to prove \eqref{eq:rev:target} for $u \in C^\infty(\Omega_R)^G$ with $\int_{\partial \Omega_R} u \, \dd A = 0$.

\medskip

\noindent \emph{Step 1. Estimates on the inner cylinder.}
Write $u(r, \theta, z) = v(r, z)$. The zero boundary mean of $u$ does not imply that its trace on $r = L$ has zero mean, so we cannot apply Lemma \ref{lem: poincare on cylinders} directly. We therefore set
\begin{align*}
	\alpha = \frac{1}{2} \int_{-1}^1 v(L, z) \, \dd z, \qquad w = u - \alpha, \qquad w(r, \theta, z) = \tilde v(r, z).
\end{align*}
Then $\int_{-1}^1 \tilde v(L, z) \, \dd z = 0$ and $\nabla w = \nabla u$. Moreover, since $\int_{\partial \Omega_R} u \, \dd A = 0$,
\begin{align}
	\int_{\partial \Omega_R} w^2 \, \dd A = \int_{\partial \Omega_R} u^2 \, \dd A + \alpha^2 \Vol(\partial \Omega_R) \geq \int_{\partial \Omega_R} u^2 \, \dd A. \label{eq:rev:shift}
\end{align}
Thus, it suffices to estimate the $L^2$ boundary norm of $w$ by its Dirichlet energy. The $G$-invariance of $w$ gives
\begin{align*}
E_R(w) = \int_{\Omega_R} \vert \nabla w \vert^2 \, \dd V = 2 \pi \int_{-1}^1 \int_0^{\lambda_R(z)} (\tilde v_r^2 + \tilde v_z^2) s \, \dd s \, \dd z.
\end{align*}

Moreover, $a < 1/2$ and $R \geq 2$ imply that $L > 1$ and $R/2 \leq L \leq R$. Using the inclusion $Q_L \subset \Omega_R$ and Lemma \ref{lem: poincare on cylinders}, we obtain
\begin{align}
	R \int_{-1}^1 \tilde v(L, z)^2 \, \dd z &\leq C E_R(w), \label{eq:rev:inner-lateral} \\
	\int_{\overline{\mathbb B_L^2} \times \{-1, 1\}} w^2 \, \dd A &\leq C R^2 E_R(w). \label{eq:rev:inner-faces}
\end{align}
Notice that the second inequality already controls the two flat portions of $\partial \Omega_R$.

\medskip

\noindent \emph{Step 2. Radial control.}
We now compare the values of $w$ in $\Omega_R \setminus Q_L$ with those on $r = L$. For each $z \in (-1, 1)$, the radial segment $[L, \lambda_R(z)]$ lies in the closure of the meridian region and has length at most $a$. 
Since $0 \leq r - L \leq \lambda_R(z) - L \leq a$, the fundamental theorem of calculus and Cauchy--Schwarz give
\begin{align*}
\vert \tilde v(r, z) \vert^2 &\leq 2 \tilde v(L, z)^2 + 2 \left\vert \int_L^r \tilde v_r(s, z) \, \dd s \right\vert^2 \\
&\leq 2 \tilde v(L, z)^2 + 2(r - L) \int_L^r \vert \tilde v_r(s, z) \vert^2 \, \dd s \\
&\leq 2 \tilde v(L, z)^2 + 2a \int_L^{\lambda_R(z)} \vert \tilde v_r(s, z) \vert^2 \, \dd s.
\end{align*}
For this fixed function $w$, define
\begin{align*}
H_R(z) = 2 \tilde v(L, z)^2 + 2a \int_L^{\lambda_R(z)} \vert \tilde v_r(s, z) \vert^2 \, \dd s.
\end{align*}
Thus $\vert \tilde v(r, z) \vert^2 \leq H_R(z)$ along the entire radial segment. We next estimate the integral of $H_R$. Since $s \geq L \geq R/2$ in the radial integral, the cylindrical expression for the Dirichlet energy gives
\begin{align*}
\int_{-1}^1 \int_L^{\lambda_R(z)} \vert \tilde v_r(s, z) \vert^2 \, \dd s \, \dd z \leq \frac{2}{R} \int_{-1}^1 \int_L^{\lambda_R(z)} \vert \tilde v_r(s, z) \vert^2 s \, \dd s \, \dd z \leq \frac{E_R(w)}{\pi R}.
\end{align*}
Recalling the definition of $H_R$ and combining the preceding estimate with \eqref{eq:rev:inner-lateral}, we obtain
\begin{align*}
\int_{-1}^1 H_R(z) \, \dd z &= 2 \int_{-1}^1 \tilde v(L, z)^2 \, \dd z + 2a \int_{-1}^1 \int_L^{\lambda_R(z)} \vert \tilde v_r(s, z) \vert^2 \, \dd s \, \dd z \\
&\leq \frac{2C}{R} E_R(w) + \frac{2a}{\pi R} E_R(w).
\end{align*}
Since $a$ is independent of $R$, we conclude that
\begin{align}
\int_{-1}^1 H_R(z) \, \dd z \leq \frac{C}{R} E_R(w). \label{eq:rev:H-bound}
\end{align}
This estimate will control the boundary values on the oscillating and vertical portions of $\partial \Omega_R$.
	
\medskip

\noindent \emph{Step 3. Estimates on the remaining boundary portions.}
We now use the preceding estimates to control the $L^2$-norm of $w$ on the oscillating, vertical and rounded portions of $\partial \Omega_R$.
	
On $\Gamma_R^{\mathrm{osc}}$, we have $\rho_R(z) = \lambda_R(z)$, so $\tilde v(\rho_R(z), z)^2 \leq H_R(z)$. Moreover, if $M_\eta = \norm{\eta'}_{L^\infty(\mathbb R)}$, then
\begin{align*}
\rho_R(z) \sqrt{1 + \rho_R'(z)^2} \leq R \sqrt{1 + R^4 M_\eta^2} \leq C R^3.
\end{align*}

Using the surface parametrisation introduced above, the $G$-invariance of $w$, and \eqref{eq:rev:H-bound}, we obtain
\begin{align}
\int_{\Gamma_R^{\mathrm{osc}}} w^2 \, \dd A &= 2 \pi \int_{-1/2}^{1/2} \tilde v(\rho_R(z), z)^2 \rho_R(z) \sqrt{1 + \rho_R'(z)^2} \, \dd z \notag \\
&\leq C R^3 \int_{-1/2}^{1/2} H_R(z) \, \dd z \leq C R^2 E_R(w). \label{eq:rev:corrugated-trace}
\end{align}

For the unchanged vertical portions $\Gamma_R^{\mathrm{vert}}$, the surface element is $R \, \dd \theta \, \dd z$. Since $\tilde v(R, z)^2 \leq H_R(z)$ there, the same argument gives
\begin{align}
\int_{\Gamma_R^{\mathrm{vert}}} w^2 \, \dd A \leq 2 \pi R \int_{-1}^1 H_R(z) \, \dd z \leq C E_R(w). \label{eq:rev:vertical-trace}
\end{align}
	
It remains to estimate the rounded portions. Consider the compact connected Lipschitz domain
\begin{align*}
D_R^+ = \left\{ (r, z) : \frac{1}{2} \leq z \leq 1, \, L \leq r \leq \lambda_R(z) \right\}
\end{align*}
and the boundary segment $\Sigma_R^+ = \{L\} \times (1/2, 1)$. Since the rounding profile is independent of $R$ up to translation in the radial direction, the pairs $(D_R^+, \Sigma_R^+)$ are all isometric by translations. Consequently, the continuity of the trace map $H^1(D_R^+) \to L^2(\gamma_R)$ and the Poincar\'e--trace inequality recalled in Section \ref{sec: Spectral convergence under conformal collapse} give
\begin{align*}
\int_{\gamma_R} \tilde v^2 \, \dd s &\leq C \norm{\tilde v}_{H^1(D_R^+)}^2 \leq C \left( \int_{D_R^+} (\tilde v_r^2 + \tilde v_z^2) \, \dd r \, \dd z + \int_{1/2}^1 \tilde v(L, z)^2 \, \dd z \right),
\end{align*}
where $\dd s$ denotes Euclidean arclength. Let $\Gamma_R^{\mathrm{round},+}$ be the surface generated by $\gamma_R$. Because $R/2 \leq r \leq R$ on $D_R^+$, integration in the angular variable yields
\begin{align*}
\int_{\Gamma_R^{\mathrm{round},+}} w^2 \, \dd A &= 2 \pi \int_{\gamma_R} \tilde v^2 r \, \dd s \\	
&\leq C R \left( \int_{D_R^+} (\tilde v_r^2 + \tilde v_z^2) \, \dd r \, \dd z + \int_{1/2}^1 \tilde v(L, z)^2 \, \dd z \right) \leq C E_R(w).
\end{align*}
Here the last inequality follows from \eqref{eq:rev:inner-lateral} and the cylindrical expression for the Dirichlet energy, since $r \geq R/2$ on $D_R^+$. Applying the same argument to the reflected lower region gives
\begin{align}
\int_{\Gamma_R^{\mathrm{round}}} w^2 \, \dd A \leq C E_R(w), \label{eq:rev:rounded-trace}
\end{align}
where $\Gamma_R^{\mathrm{round}}$ denotes the union of the two rounded portions.
	
Combining \eqref{eq:rev:inner-faces}, \eqref{eq:rev:corrugated-trace}, \eqref{eq:rev:vertical-trace} and \eqref{eq:rev:rounded-trace}, we obtain
\begin{align*}
\int_{\partial \Omega_R} w^2 \, \dd A \leq C R^2 E_R(w).
\end{align*}
Together with \eqref{eq:rev:shift} and $\nabla w = \nabla u$, this proves \eqref{eq:rev:target} for smooth $G$-invariant functions with zero boundary mean, and hence for all such functions in $H^1(\Omega_R)^G$ by density. The variational characterisation of the first positive $G$-invariant Steklov eigenvalue then yields
\begin{align}
\sigma_1^G(\Omega_R) \geq \frac{c}{R^2}. \label{eq:rev:sigma-bound}
\end{align}

Finally, combining the volume bound \eqref{eq:rev:volume-bound}, the boundary area bound \eqref{eq:rev:area-bound}, and the eigenvalue bound \eqref{eq:rev:sigma-bound}, we obtain
\begin{align*}
\overline{\sigma}_1^G(\Omega_R) = \sigma_1^G(\Omega_R) \Vol(\Omega_R)^{-1/3} \Vol(\partial \Omega_R) \geq c R^{1/3} \longrightarrow +\infty
\end{align*}
as $R \to +\infty$.	This proves the theorem in dimension three.

\medskip

\noindent \emph{Higher dimensions.}
Let $n \geq 4$ and let $G = \mathrm{O}(n - 1)$ act on $\mathbb R^n = \mathbb R^{n - 1} \times \mathbb R$ by orthogonal transformations on the first factor. For $s > 0$, let $\mathbb B_s^{n - 1}$ denote the Euclidean ball of radius $s$ in $\mathbb R^{n - 1}$, and set $Q_s = \overline{\mathbb B_s^{n - 1}} \times [-1, 1]$. The same meridian profiles define compact smooth $G$-invariant domains $\Omega_R$ diffeomorphic to a closed ball, with $Q_L \subset \Omega_R \subset Q_R$.

The inclusion $\Omega_R \subset Q_R$ immediately gives
\begin{align*}
\Vol(\Omega_R) \leq C R^{n - 1}.
\end{align*}
Moreover, integrating the area element of the oscillating boundary and using the same periodicity argument as in dimension three, we obtain
\begin{align*}
\Vol(\partial \Omega_R) \geq c R^{n - 2} \int_{-1/2}^{1/2} \vert \rho_R'(z) \vert \, \dd z \geq c R^n.
\end{align*}

It remains to estimate $\sigma_1^G(\Omega_R)$. Let $u \in C^\infty(\Omega_R)^G$ have zero boundary mean, and write $u(x,z) = v(r,z)$, where $r = \vert x \vert$. As in Step 1, set
\begin{align*}
\alpha = \frac{1}{2} \int_{-1}^1 v(L,z) \, \dd z, \qquad w = u - \alpha, \qquad w(x,z) = \tilde v(r,z),
\end{align*}
and let $E_R(w) = \int_{\Omega_R} \vert \nabla w \vert^2 \, \dd V$. Then $\int_{-1}^1 \tilde v(L,z) \, \dd z = 0$, $\nabla w = \nabla u$, and, as in Step 1,
\begin{align*}
\int_{\partial \Omega_R} u^2 \, \dd A \leq \int_{\partial \Omega_R} w^2 \, \dd A.
\end{align*}
Define $H_R$ as in Step 2, namely
\begin{align*}
H_R(z) = 2 \tilde v(L,z)^2 + 2a \int_L^{\lambda_R(z)} \vert \tilde v_r(s,z) \vert^2 \, \dd s.
\end{align*}
The proof of Lemma \ref{lem: poincare on cylinders} carries over with radial weight $r^{n - 2}$, giving
\begin{align*}
R^{n - 2} \int_{-1}^1 \tilde v(L, z)^2 \, \dd z \leq C E_R(w), \qquad \int_{\overline{\mathbb B_L^{n - 1}} \times \{-1, 1\}} w^2 \, \dd A \leq C R^2 E_R(w).
\end{align*}
The argument of Step 2 then yields $\int_{-1}^1 H_R(z) \, \dd z \leq C R^{2 - n} E_R(w)$. 
Since the surface density on the oscillating portion is bounded by $C R^n$, the estimate for $H_R$ gives
\begin{align*}
\int_{\Gamma_R^{\mathrm{osc}}} w^2 \, \dd A \leq C R^2 E_R(w).
\end{align*}
The vertical and rounded portions contribute at most $C E_R(w)$, by the same arguments as in dimension three, using that the radial weight is comparable to $R^{n - 2}$ in the outer region. Consequently,
\begin{align*}
\int_{\partial \Omega_R} u^2 \, \dd A \leq C R^2 \int_{\Omega_R} \vert \nabla u \vert^2 \, \dd V,
\end{align*}
and hence 
\begin{align*}
\sigma_1^G(\Omega_R) \geq \frac{c}{R^2}.
\end{align*}

Finally, combining the three bounds, we obtain
\begin{align*}
\overline{\sigma}_1^G(\Omega_R) = \sigma_1^G(\Omega_R) \Vol(\Omega_R)^{\frac{2 - n}{n}} \Vol(\partial \Omega_R) \geq c R^{\frac{n - 2}{n}} \longrightarrow +\infty
\end{align*}
as $R \to +\infty$. This completes the proof.
\end{proof}


\begin{thebibliography}{99}
\bibitem{abreu}
M.  Abreu and P. Freitas,
\emph{On the invariant spectrum of $S^1$-invariant metrics on $S^2$},
Proc. London Math. Soc. (3) \textbf{84} (1) (2002), 213--230.

\bibitem{colbois survey}
B. Colbois, A. Girouard, C. Gordon, and D. Sher,
\emph{Some recent developments on the Steklov eigenvalue problem},
Rev. Mat. Complut. \textbf{37} (2024), 1--161.

\bibitem{colbois_G}
B. Colbois, E. Dryden, and A. El Soufi, 
\emph{Extremal $G$-invariant eigenvalues of the Laplacian of $G$-invariant metrics},
Math. Z. \textbf{258} (2008), 29--41.

\bibitem{colin_verdiere}
Y. Colin de Verdi\`ere,
\emph{Sur la multiplicit\'e de la premi\`ere valeur propre non nulle du Laplacien},
Comment. Math. Helv. \textbf{61} (1986), 254--270.

\bibitem{enciso-pistoia-provenzano}
A. Enciso, A. Pistoia and L. Provenzano,
\emph{Nonexistence of Courant-type nodal domain bounds for eigenfunctions of the Dirichlet-to-Neumann operator},
J. Eur. Math. Soc. (2025), published online first.

\bibitem{engman}
M. Engman,
\emph{Trace formulae for $S^1$ invariant Green's operators on $S^2$},
Manuscripta Math. \textbf{93} (1997), 357--368.

\bibitem{fraser2011}
A. Fraser and R. Schoen,
\emph{The first Steklov eigenvalue, conformal geometry, and minimal surfaces},
Adv. Math. \textbf{226} (5) (2011), 4011--4030.

\bibitem{fraser}
A. Fraser and R. Schoen,
\emph{Minimal surfaces and eigenvalue problems},
Contemp. Math. \textbf{599} (2013), 105--121.

\bibitem{girouard2010}
A. Girouard and I. Polterovich,
\emph{On the Hersch--Payne--Schiffer inequalities for Steklov eigenvalues},
Funktsional. Anal. i Prilozhen. \textbf{44} (2) (2010), 33--47.

\bibitem{girouard2012}
A. Girouard and I. Polterovich,
\emph{Upper bounds for Steklov eigenvalues on surfaces},
Electron. Res. Announc. Math. Sci. \textbf{19} (2012), 77--85.

\bibitem{girouard}
A. Girouard and I. Polterovich,
\emph{Spectral geometry of the Steklov problem},
J. Spectr. Theory \textbf{7} (2) (2017), 321--359.

\bibitem{girouard2021}
A. Girouard, M. Karpukhin, and J. Lagac\'e,
\emph{Continuity of eigenvalues and shape optimisation for Laplace and Steklov problems},
Geom. Funct. Anal. \textbf{31} (2021), 513--561.

\bibitem{hassannezhad}
A. Hassannezhad,
\emph{Conformal upper bounds for the eigenvalues of the Laplacian and Steklov problem},
J. Funct. Anal. \textbf{261} (2011), 3419--3436.

\bibitem{hersch}
J. Hersch, L. E. Payne, and M. M. Schiffer,
\emph{Some inequalities for Stekloff eigenvalues},
Arch. Ration. Mech. Anal. \textbf{57} (1974), 99--114.

\bibitem{karpukhin2017}
M. Karpukhin,
\emph{Bounds between Laplace and Steklov eigenvalues on nonnegatively curved manifolds},
Electron. Res. Announc. Math. Sci. \textbf{24} (2017), 100--109.

\bibitem{karpukhin2022}
M. Karpukhin and A. M\'etras,
\emph{Laplace and Steklov Extremal Metrics via $n$-Harmonic Maps},
J. Geom. Anal. \textbf{32}, 154 (2022). 

\bibitem{kriegl}
A. Kriegl, P. Michor and A. Rainer,
\emph{Many parameter H\"older perturbation of unbounded operators},
Math. Ann. \textbf{353} (2012), 519--522.

\bibitem{korevaar}
N. Korevaar,
\emph{Upper bounds for eigenvalues of conformal metrics},
J. Differential Geom. \textbf{37}(1) (1993), 73--93. 

\bibitem{kuznetsov}
N. Kuznetsov, T. Kulczycki, M. Kwa\'snicki, A. Nazarov, S. Poborchi, I. Polterovich, and B. Siudeja,
\emph{The legacy of Vladimir Andreevich Steklov},
Notices Am. Math. Soc. \textbf{61} (1) (2014), 9--22.

\bibitem{nadirashvili}
N. Nadirashvili,
\emph{Berger's isoperimetric problem and minimal immersions of surfaces},
Geom. Funct. Anal. \textbf{6} (5) (1996), 877--897.

\bibitem{weinstock}
R. Weinstock,
\emph{Inequalities for a classical eigenvalue problem},
J. Rational Mech. Anal. \textbf{3} (1954), 745--753.
\end{thebibliography}
\end{document}